\documentclass[11pt]{amsart}
\usepackage{amsmath,amssymb,amsthm}
\usepackage{microtype}
\usepackage{enumitem}
\usepackage[svgnames]{xcolor}
\usepackage[colorlinks=true,linkcolor=Navy,citecolor=Navy,urlcolor=Navy]{hyperref}
\usepackage[margin=1.08in]{geometry}

\newtheorem{theorem}{Theorem}[section]
\newtheorem{proposition}[theorem]{Proposition}
\newtheorem{lemma}[theorem]{Lemma}
\newtheorem{corollary}[theorem]{Corollary}
\newtheorem{remark}[theorem]{Remark}
\newtheorem{definition}[theorem]{Definition}
\newtheorem{maintheorem}{Theorem}

\title[Permutation polynomials of given index]{Permutation polynomials of given index up to quasi-multiplicative equivalence}
\author{Xiang Fan}
\address{School of Mathematics, Sun Yat-sen University, Guangzhou 510275, China}
\email{fanx8@mail.sysu.edu.cn}
\hypersetup{
  pdftitle={Permutation polynomials of given index up to quasi-multiplicative equivalence},
  pdfauthor={Xiang Fan},
  pdfsubject={Classification and enumeration of permutation polynomials of given index up to quasi-multiplicative equivalence},
  pdfkeywords={permutation polynomial, index of a polynomial, quasi-multiplicative equivalence, cyclotomic mapping, double coset, quasi-polynomial, Jacobi sum, elliptic curve over a finite field}
}

\date{November 2026}

\begin{document}
\begin{abstract}
Let $q$ be a prime power and let $\ell$ be a divisor of $q-1$.  We classify the permutation polynomials of $\mathbb F_q$ of polynomial index exactly $\ell$ up to quasi-multiplicative equivalence $G(X)\equiv uF(vX^d)\pmod{X^q-X}$, where $u,v\in\mathbb F_q^*$ and $\gcd(d,q-1)=1$, and we count the classes.  For polynomials vanishing at $0$, the classes are the orbits of an explicit finite group on vectors of values at the $\ell$-th roots of unity.  Through the permutation it induces on the roots of unity, each orbit lies over a double coset of the symmetric group $S_\ell$ with respect to its translation subgroup and an explicit affine group, and the stabilizers are explicit.  For polynomials with nonzero constant term the action is free.  Burnside's lemma, M\"obius inversion over the coefficient support and congruences along the cycles of affine permutations give an exact formula for the number of classes.  For fixed $\ell\geqslant2$ these numbers are quasi-polynomials in $s=(q-1)/\ell$ of degree at most $\ell-1$, with period dividing $\ell$ for odd $\ell$ and $2\ell$ for even $\ell$; for polynomials vanishing at $0$ the degree is exactly $\ell-1$, and we determine the leading coefficient.  For odd prime $\ell$ and for $\ell=4$ the counts are given in closed form.  For index four we also count the binomials and trinomials at each admissible coefficient support.  The sums of two binomial counts that occur are expressed through quartic Jacobi sums and Frobenius traces of elliptic curves and of a curve of genus two; one of the elliptic curves is the reduction of the curve $Y^2=X^3+8X^2+32X$, which has no complex multiplication over $\mathbb Q$.  The trinomial counts are expressed through quartic cyclotomic numbers and these sums.  All these terms cancel when the counts with a fixed residue of the vanishing order modulo $s$ are added.
\end{abstract}
\maketitle

\section{Introduction}\label{sec:intro}

\subsection{Permutation polynomials of given index}\label{subsec:index}
Let $q$ be a prime power.  A polynomial $F\in\mathbb F_q[X]$ is a \emph{permutation polynomial} (PP) of $\mathbb F_q$ if the map $x\mapsto F(x)$ is a bijection of $\mathbb F_q$.  Since $X^q-X$ vanishes on $\mathbb F_q$, every polynomial function on $\mathbb F_q$ has a unique representative of degree at most $q-1$.  We always work with this reduced representative, and we identify polynomials that induce the same function on $\mathbb F_q$.  A nonconstant reduced polynomial is written as
\[
 F(X)=\beta+\sum_{i=1}^{n}c_iX^{e_i},\qquad
 1\leqslant e_1<\cdots<e_n\leqslant q-1,\quad c_i\neq0,
\]
and $r:=e_1$ denotes the vanishing order of $F-\beta$ at $0$.

\begin{definition}\label{def:index}
The \emph{polynomial index} $\operatorname{ind}(F)$ of $F$ is the least divisor $\ell$ of $q-1$ for which $e_1,\dots,e_n$ are congruent modulo $(q-1)/\ell$.
\end{definition}

The index can be read off from the exponents.  Put
\[
 \Delta(F):=\gcd(q-1,\,e_2-r,\dots,e_n-r),
\]
with $\Delta(F):=q-1$ when $n=1$.  Then $\Delta(F)$ is the largest divisor of $q-1$ that divides every exponent difference, so $\operatorname{ind}(F)=(q-1)/\Delta(F)$.  In particular, the index is determined by the polynomial function $F$ and not by a choice of presentation, and a monomial has index $1$.

Let $\ell$ be a divisor of $q-1$ and put $s:=(q-1)/\ell$.  A presentation of $F$ in terms of $X^s$ exists precisely when $\operatorname{ind}(F)\mid\ell$.  In that case there are unique polynomials $H,h$ and a unique integer $a\in\{1,\dots,s\}$ such that
\[
 F(X)=\beta+X^rH(X^s)=\beta+X^ah(X^s),
 \qquad H(0)\neq0,\quad a\equiv r\pmod s,\quad \deg H,\ \deg h<\ell,
\]
namely $h(T)=T^{(r-a)/s}H(T)$.  The first presentation, in terms of the vanishing order $r$, is the one used in \cite[p.~196]{AGW2009} and \cite[\S15.1, pp.~319--320]{Wang2019}.  The second, in terms of the \emph{residue exponent} $a$, is the one adapted to the normalization used below; note that $h(0)$ may vanish.  Whether $F$ permutes $\mathbb F_q$ is decided by the condition $\gcd(a,s)=1$ together with a map of the group $\mu_\ell$ of $\ell$-th roots of unity to itself \cite[Lemma~2.1 and Theorem~2.2]{AkbaryWang2007}; see also \cite{WanLidl1991} and Lemma~\ref{lem:cyclo} below.

Li, Guan and Wu \cite{LiGuanWu2025} determined the number of PPs of $\mathbb F_q$ of exact index $\ell$.  Their count does not identify polynomials related by the evident symmetries of the problem: replacing $X$ by $vX^d$ with $v\in\mathbb F_q^*$ and $\gcd(d,q-1)=1$, or multiplying by a nonzero constant, preserves both the permutation property and the index.  We study PPs of a given index up to these symmetries.

\begin{definition}[Quasi-multiplicative equivalence]\label{def:qm-equivalence}
Two polynomials $F,G\in\mathbb F_q[X]$ are said to be \emph{quasi-multiplicatively equivalent}, or \emph{QM equivalent}, written $F\sim_{\mathrm{QM}}G$, if there are $u,v\in\mathbb F_q^*$ and a positive integer $d$ with $\gcd(d,q-1)=1$ such that
\[
 G(X)\equiv uF(vX^d)\pmod{X^q-X}.
\]
\end{definition}

Equivalently, $G(x)=uF(vx^d)$ for every $x\in\mathbb F_q$; thus QM equivalence is a relation between polynomial functions.  The substitutions $x\mapsto vx^d$ are closed under composition and inversion, and so are the output scalings; hence $\sim_{\mathrm{QM}}$ is an equivalence relation, and it preserves the permutation property.  On reduced representatives, the substitution $X\mapsto vX^d$ multiplies the nonconstant exponents by the unit $d$ modulo $q-1$; it therefore leaves $\Delta(F)$, and with it the index, unchanged.  Finally $G(0)=uF(0)$, so the conditions $F(0)=0$ and $F(0)\neq0$ are invariant under QM equivalence.

Transformations of this type are called quasi-multiplicative equivalences by Wu, Yuan, Ding and Ma \cite[Definition~2.2]{WuYuanDingMa2017}, while Ding and Zieve use the term multiplicative equivalence for the relation modulo $X^q-X$ \cite[Definition~1.4 and Remark~1.5]{DingZieve2023}.  As the latter term has also been used for narrower notions, we keep the name QM equivalence.

\subsection{Main results}\label{subsec:main}
Throughout, $\ell$ is a divisor of $q-1$,
\[
 s=\frac{q-1}{\ell},\qquad g:=\gcd(\ell,s),
\]
$\gamma$ is a fixed generator of $\mathbb F_q^*$, $\zeta:=\gamma^s$, and $\mu_\ell:=\langle\zeta\rangle$.  Let $S_\ell$ be the symmetric group of $\mathbb Z/\ell\mathbb Z$, let $\mathcal T_\ell$ be its subgroup of translations, and put
\[
 H_{\ell,s}:=\{\delta\in(\mathbb Z/\ell\mathbb Z)^*:\delta\equiv1\pmod g\},\qquad
 \mathcal A_{\ell,s}:=(\mathbb Z/\ell\mathbb Z)\rtimes H_{\ell,s};
\]
thus $\mathcal A_{\ell,s}$ is the group of affine permutations $j\mapsto t+\delta j$ of $\mathbb Z/\ell\mathbb Z$ with $t\in\mathbb Z/\ell\mathbb Z$ and $\delta\in H_{\ell,s}$.  If $F=\beta+X^ah(X^s)$ is a PP, the map $u\mapsto u^ah(u)^s$ permutes $\mu_\ell$ (Lemma~\ref{lem:cyclo}), so there is a unique $\sigma_F\in S_\ell$ with
\[
 \zeta^{aj}h(\zeta^j)^s=\zeta^{\sigma_F(j)}\qquad(j\in\mathbb Z/\ell\mathbb Z).
\]
We call $\sigma_F$ the \emph{induced permutation} of $F$.  The first theorem describes the QM classes.

\begin{maintheorem}\label{thm:A}
Let $\ell\mid q-1$.
\begin{enumerate}[label=(\arabic*)]
\item\label{A:normal} Every QM class of PPs of exact index $\ell$ vanishing at $0$ contains a polynomial $Xh(X^s)$ with $h(1)=1$.  On these normalized polynomials, QM equivalence is the orbit relation of an explicit group $\Gamma^0_{\ell,s}$ of order $\ell g\varphi(\ell)/\varphi(g)$ acting on the value vectors $(h(\zeta^j))_{j\in\mathbb Z/\ell\mathbb Z}$.
\item\label{A:double} There is a surjective homomorphism $\pi\colon\Gamma^0_{\ell,s}\to\mathcal A_{\ell,s}$ with kernel of order $g$ such that, if $\omega\in\Gamma^0_{\ell,s}$ maps the value vector of a normalized polynomial $F$ to that of $F'$ and $\pi(\omega)$ is $j\mapsto t+\delta j$, then $\sigma_{F'}(j)=\sigma_F(t+\delta j)-\sigma_F(t)$.  Consequently, through the induced permutations of its normalized members, each class in \ref{A:normal} lies over a double coset in $\mathcal T_\ell\backslash S_\ell/\mathcal A_{\ell,s}$, and those lying over the double coset of a permutation $\sigma$ with $\sigma(0)=0$ are the orbits of
\[
 K_\sigma:=\pi^{-1}\bigl(\mathcal A_{\ell,s}\cap\sigma^{-1}\mathcal T_\ell\sigma\bigr)
\]
on the normalized polynomials with induced permutation $\sigma$.  A double coset may carry no polynomial of exact index $\ell$.
\item\label{A:fixed} If all cycles of $\alpha\in\mathcal A_{\ell,s}$ have a common length $L$, then right multiplication by $\alpha$ fixes exactly $\varphi(L)L^{m-1}(m-1)!$ left cosets $\mathcal T_\ell\sigma$, where $m=\ell/L$; otherwise it fixes none.
\item\label{A:nonzero} Every QM class of PPs of exact index $\ell$ with $F(0)\neq0$ contains polynomials with constant term $1$.  Two such polynomials are QM equivalent if and only if they differ by a substitution $X\mapsto vX^d$, and the group of these substitutions acts freely on the PPs with constant term $1$.
\end{enumerate}
In both cases, each class has a unique lexicographically least normalized representative in explicit coordinates \textup{(}Theorems~\ref{thm:zero-quotient} and~\ref{thm:nonzero-structure}\textup{)}.
\end{maintheorem}

The group $\Gamma^0_{\ell,s}$ and its action are given in Section~\ref{sec:orbits}.  The normalized polynomials with a prescribed induced permutation are parametrized by their values, and exactness of the index is a condition on the inverse discrete Fourier transform of these values (Theorem~\ref{thm:fourier-exact}).  For counting, we write $\mathfrak Q_\ell(q)$ and $\mathfrak R_\ell(q)$ for the numbers of QM classes of PPs of $\mathbb F_q$ of exact index $\ell$ with $F(0)=0$ and with $F(0)\neq0$, respectively, and put $\mathfrak A_\ell(q):=\mathfrak Q_\ell(q)+\mathfrak R_\ell(q)$.  A \emph{quasi-polynomial} in $s$ is a function that is polynomial on each residue class modulo a fixed period.

\begin{maintheorem}\label{thm:B}
Let $\ell\mid q-1$.
\begin{enumerate}[label=(\arabic*)]
\item\label{B:formula} $\mathfrak Q_\ell(q)$ is given by an explicit finite formula \textup{(}Theorem~\ref{thm:class-number}\textup{)}: an average over the stabilizers $K_\sigma$ of M\"obius sums of the numbers of solutions of systems of linear congruences modulo $s$.  Each such number is either $0$ or $s^{N-1}\gcd(\lambda,s)$, where $N$ and $\lambda$ are the number of cycles and the greatest common divisor of the cycle lengths of an explicit affine permutation.
\item\label{B:nonzero} $\mathfrak R_\ell(q)=N_0(\ell,q)/\bigl((q-1)\varphi(q-1)\bigr)$, where
\[
 N_0(\ell,q)=\sum_{d\mid\ell}\mu(\ell/d)\,d!\left(\frac{q-1}{d}\right)^d\varphi\!\left(\frac{q-1}{d}\right)
\]
is the number of PPs of exact index $\ell$ vanishing at $0$, determined by Li, Guan and Wu \cite{LiGuanWu2025}.
\item\label{B:quasi} Fix $\ell\geqslant2$, and let $P_\ell:=\ell$ if $\ell$ is odd and $P_\ell:=2\ell$ if $\ell$ is even.  For each $\varrho\in\mathbb Z/P_\ell\mathbb Z$ there are polynomials $Q_{\ell,\varrho},R_{\ell,\varrho}\in\mathbb Q[X]$ of degree at most $\ell-1$ such that
\[
 \mathfrak Q_\ell(q)=Q_{\ell,\varrho}(s),\qquad \mathfrak R_\ell(q)=R_{\ell,\varrho}(s)
\]
whenever $q=\ell s+1$ is a prime power with $s\equiv\varrho\pmod{P_\ell}$; moreover $R_{\ell,\varrho}$ depends only on $\varrho$ modulo $\operatorname{rad}(\ell)$.  The coefficient of $X^{\ell-1}$ in $Q_{\ell,\varrho}$ is
\[
 \frac{(\ell-1)!\,\varphi(g)}{\ell g\,\varphi(\ell)},\qquad g=\gcd(\ell,\varrho).
\]
\end{enumerate}
\end{maintheorem}

Thus $\mathfrak Q_\ell$, $\mathfrak R_\ell$ and $\mathfrak A_\ell$ are quasi-polynomials in $s$ of degree at most $\ell-1$ whose periods divide $P_\ell$, and $\mathfrak Q_\ell$ has degree exactly $\ell-1$.  For odd prime indices and for index four, the formula of Theorem~\ref{thm:B}~\ref{B:formula} can be evaluated in closed form.

\begin{maintheorem}\label{thm:C}
\begin{enumerate}[label=(\arabic*)]
\item\label{C:prime} Let $p$ be an odd prime and $q-1=ps$.  Then
\[
 \mathfrak Q_p(q)=
 \begin{cases}
 \dfrac{(p-2)!\,s^{p-1}-1}{p},&p\nmid s,\\[2mm]
 \dfrac{(p-1)!\,s^{p-1}}{p^2}+p-2,&p\mid s,
 \end{cases}
 \qquad
 \mathfrak R_p(q)=
 \begin{cases}
 (p-2)!\,s^{p-1}-1,&p\nmid s,\\[1mm]
 \dfrac{(p-1)!\,s^{p-1}}{p}-1,&p\mid s.
 \end{cases}
\]
\item\label{C:four} If $q-1=4s$, then
\[
 \mathfrak Q_4(q)=
 \begin{cases}
 (3s^3+s)/4,&s\text{ odd},\\
 (3s^3+4s)/8,&s\equiv2\pmod4,\\
 3s^3/8+1,&4\mid s,
 \end{cases}
 \qquad
 \mathfrak A_4(q)=
 \begin{cases}
 (15s^3-3s)/4,&s\text{ odd},\\
 (15s^3-4s)/8,&s\equiv2\pmod4,\\
 (15s^3-8s+8)/8,&4\mid s.
 \end{cases}
\]
\end{enumerate}
\end{maintheorem}

For index four, the polynomials can also be counted according to their coefficient support.  Let $q\equiv1\pmod4$ and $s=(q-1)/4$.  For an integer $t$, let $B(t)$ be the number of $b\in\mathbb F_q^*$ such that $u\mapsto u^t(u+b)^s$ permutes $\mu_4$; here $u^t$ is defined for every integer $t$ because $\mu_4\subset\mathbb F_q^*$.  By Lemma~\ref{lem:cyclo}, if $t\geqslant1$ and $\gcd(t,s)=1$, then $B(t)$ is the number of $b\in\mathbb F_q^*$ for which $X^t(X^s+b)$ is a PP.  For $1\leqslant t\leqslant2s$ with $\gcd(t,s)=1$, let $\mathfrak T(t)$ be the number of points $[c_0:c_1:c_2]\in\mathbb P^2(\mathbb F_q)$ with $c_0c_1c_2\neq0$ such that $X^t(c_0+c_1X^s+c_2X^{2s})$ is a PP.  Let $\psi$ be the quartic character of $\mathbb F_q^*$ with $\psi(\gamma)=i$, extended by $\psi(0):=0$, put $J:=\sum_{x\in\mathbb F_q}\psi(x)\psi(1-x)$, and for a smooth projective curve $C$ over $\mathbb F_q$ put $a_q(C):=q+1-\lvert C(\mathbb F_q)\rvert$.  We use the elliptic curves
\[
 E_{128}\colon Y^2=X^3+8X^2+32X,\qquad
 E_{1728}\colon Y^2=X^3+4X,\qquad
 E_{8000}\colon Y^2=X^3+2X^2-2X-4,
\]
whose $j$-invariants are $128$, $1728$ and $8000$.  These equations have integer coefficients and define elliptic curves over $\mathbb Q$ and, by reduction, over $\mathbb F_q$; the traces $a_q(E_{128})$, $a_q(E_{1728})$ and $a_q(E_{8000})$ refer to the curves over $\mathbb F_q$.  Over $\mathbb Q$, the curves $E_{1728}$ and $E_{8000}$ have complex multiplication, whereas $E_{128}$ has none (Section~\ref{subsec:support4}).

\begin{maintheorem}\label{thm:D}
Let $q\equiv1\pmod4$, $s=(q-1)/4$, and let $r$ be an integer with $\gcd(r,s)=1$.
\begin{enumerate}[label=(\arabic*)]
\item\label{D:even} If $q\equiv1\pmod8$, then $\varepsilon:=i^s$ and $\rho:=\psi(2)$ lie in $\{\pm1\}$, and
\[
 B(r)+B(r+s)=\frac{3q-69}{16}+\frac{\rho}{2}(\operatorname{Re}J-1)-\frac{\varepsilon(\rho+1)}2
 -\frac14a_q(E_{128})-\frac9{16}a_q(E_{1728})-\frac{2-\varepsilon}{8}a_q(E_{8000}).
\]
\item\label{D:odd} If $q\equiv5\pmod8$, then $B(r)+B(r+s)=M_q+\tfrac12\Lambda_q\operatorname{Re}(i^{r+s})$, where $M_q$ is an explicit linear combination of $1$, $q$, $\operatorname{Re}\bigl(\psi(2)(J+1)\bigr)$, $a_q(E_{128})$ and $a_q(E_{1728})$, and $\Lambda_q$ is an integer expressed through $J$ and the Frobenius trace of an explicit curve of genus two \textup{(}Theorem~\ref{thm:odd-boundary}\textup{)}.
\item\label{D:tri} If $1\leqslant r\leqslant2s$, then
\[
 \mathfrak T(r)=D_r+4m_+m_-+\delta_r(q-1)-B(r)-B(r+s)-\iota_r,
\]
where $D_r$ and $m_\pm$ are explicit expressions in the quartic cyclotomic numbers of $\mathbb F_q$ \textup{(}Section~\ref{subsec:opposite}\textup{)}, $\delta_r=1$ if $r\equiv3s-2\pmod4$ and $\delta_r=0$ otherwise, and $\iota_r=1$ if $r+s$ is odd and $\iota_r=0$ otherwise \textup{(}Theorem~\ref{thm:trinomial-count}\textup{)}.
\item\label{D:strata} Let $1\leqslant a\leqslant s$ with $\gcd(a,s)=1$, and put $L_s:=s(3+(-1)^s)$.  Counted up to multiplication by nonzero constants, the PPs of exact index four vanishing at $0$ with vanishing order $a$, $a+s$ and $a+2s$ number
\[
 6s^3-L_s-\mathfrak T(a+s)-B(a+s)-B(a+2s),\qquad \mathfrak T(a+s)+B(a+s),\qquad B(a+2s),
\]
respectively, and none has vanishing order $a+3s$.  In particular, the total over the three vanishing orders is $6s^3-L_s$.
\end{enumerate}
\end{maintheorem}

By parts \ref{D:even} and~\ref{D:odd}, the sums $B(r)+B(r+s)$ are expressed through the Jacobi sum $J$, the Frobenius traces of the elliptic curves above and, when $q\equiv5\pmod8$, the Frobenius trace of a curve of genus two; by part \ref{D:tri}, the trinomial counts $\mathfrak T(r)$ are expressed through quartic cyclotomic numbers and these sums.  Each individual number $B(t)$ is a finite character sum, see \eqref{eq:B-fourier}.  These sums determine the individual numbers $B(t)$ unless $q\equiv9\pmod{16}$; in that case only $B(t)+B(t+s)$ is evaluated, and the counts in part \ref{D:strata} still involve individual binomial counts (Remark~\ref{rem:individual-B}).  Part \ref{D:strata} shows that all these terms disappear from the total over a fixed residue exponent.  In particular, the Frobenius trace of $E_{128}$, a curve without complex multiplication over $\mathbb Q$, appears in the formulas for $B(r)+B(r+s)$ and $\mathfrak T(r)$ but not in their total.  Summing part \ref{D:strata} over the residue exponents recovers the count $N_0(4,q)$ of Li, Guan and Wu (Corollary~\ref{cor:LGW-recovery}).

\subsection{Outline of the proofs}\label{subsec:method}
Lemma~\ref{lem:cyclo} reduces the permutation property of $X^ah(X^s)$ to a condition on the values $A_j=h(\zeta^j)$.  The coefficients of $h$ are recovered from these values by the inverse discrete Fourier transform, and the index is exactly $\ell$ if and only if there are no prime $p\mid\ell$ and no $b$ for which the eigenvalue relations $A_{j+\ell/p}=\zeta^{b\ell/p}A_j$ hold for all $j$ (Theorem~\ref{thm:fourier-exact}).  Normalizing the residue exponent to $1$ and the value at $1$ to $1$ leaves a finite group acting on the remaining values, and its action on induced permutations is affine (Lemma~\ref{lem:qm-action}); this gives Theorem~\ref{thm:A}.

The counting problem is more delicate, because exactness of the index is a property of the coefficient vector and not of the induced permutation.  We impose it by M\"obius inversion over the divisors of $\ell$.  In every term of the resulting sum over fixed points, the fixed-point equations together with the support conditions become a system of linear congruences modulo $s$ along the cycles of an affine permutation (Lemma~\ref{lem:cycle-count}).  The solution set is either empty or has $s^{N-1}\gcd(\lambda,s)$ elements, and its nonemptiness depends on $s$ only through a residue class modulo $P_\ell$; the latter rests on a divisibility property of the cycle sums of affine permutations (Lemma~\ref{lem:affine-cycle-sum}).  This yields Theorem~\ref{thm:B}.  For odd prime indices the remaining stabilizer terms are controlled by quadratic Gauss sums, which gives Theorem~\ref{thm:C}~\ref{C:prime}.

For index four the global count is elementary, but the counts with a prescribed coefficient support are governed by character sums.  For binomials we expand the condition that the four values $u^t(u+b)^s$, $u\in\mu_4$, form a permutation of $\mu_4$ by Fourier analysis on $(\mathbb Z/4\mathbb Z)^4$.  Only $64$ Fourier coefficients survive, and the resulting character sums are identified with Frobenius traces on Kummer sheaves and then, through quotients of cyclic covers by involutions, with Frobenius traces of the curves above (Section~\ref{sec:binomial}).  For trinomials, a sum--product characterization of $\mu_4$ and a decomposition according to opposite pairs in $\mu_4$ reduce the count to quartic cyclotomic numbers and to the binomial counts (Section~\ref{sec:trinomial}).

\subsection{Relation to previous work}\label{subsec:literature}
Presentations $X^rH(X^{(q-1)/\ell})$ and the reduction of the permutation property to $\mu_\ell$ are studied in \cite{AGW2009,AkbaryWang2007,WanLidl1991}; see \cite{Wang2019} for an account of polynomials over finite fields organized around the index.  Li, Guan and Wu \cite{LiGuanWu2025} enumerate the PPs of exact index $\ell$ before passing to QM classes; we use their count for polynomials with nonzero constant term (Theorem~\ref{thm:B}~\ref{B:nonzero}) and recover it for index four from the counts by vanishing order (Corollary~\ref{cor:LGW-recovery}).

Bors and Wang \cite[Theorem~1.1 and \S6]{BorsWang2022} describe generalized cyclotomic permutation groups through wreath products and study their cycle structure.  A chosen cyclotomic partition index need not coincide with the least polynomial index used here, and conjugacy in those groups is different from QM equivalence.  The double-coset argument of Section~\ref{sec:orbits} uses only standard group theory; what is specific to the present setting is its combination with the exact-index condition on the value vectors over each induced permutation and with the fixed-point counts of the stabilizers.

Hou and Pallozzi Lavorante \cite[\S2, Theorem~2.2]{HouPallozziLavorante2023} give canonical forms for permutation binomials under an equivalence that also allows the output Frobenius transformation $F\mapsto F^p$, where $p$ is the characteristic.  For a reduced binomial $aX^m+bX^n$, their parameter $D=\gcd(q-1,m-n)$ corresponds to the polynomial index $(q-1)/D$.  Their canonical forms concern the enlarged equivalence and do not directly classify the QM orbits considered here.  Frameworks for counting equivalence classes of polynomials and of rational functions appear in \cite{Hou2025,Mullen1976}.  Yadav, Gupta, Singh and Yadav give an algorithm that decides whether two permutation trinomials are QM equivalent and apply it to specified families \cite{YadavGuptaSinghYadav2026}; QM equivalence within specified families is also examined in \cite{YadavSinghGupta2025}.  Lee and Park give criteria for PPs of index three \cite{LeePark1997}, and Park and Lee give a criterion for PPs of the form $X^rh(X^s)$ in terms of the subgroup of $s$-th powers \cite[Theorem~2.3]{ParkLee2001}; constructions of PPs of index $q+1$ over $\mathbb F_{q^2}$ appear in \cite{FengWang2026}.  To the best of our knowledge, none of these works gives the complete classification of the PPs of exact index $\ell$ up to QM equivalence for arbitrary $\ell$.

For the arithmetic of index four, Oliveira and Brochero Mart\'inez \cite{OliveiraBrochero2022} obtained exact counts of permutation binomials of indices $2$ and $3$ by counting points on curves, with an elliptic-curve calculation in the latter case.

\subsection{Organization}\label{subsec:organization}
Section~\ref{sec:prelim} translates the index and the permutation property into conditions on values at roots of unity.  Section~\ref{sec:orbits} proves Theorem~\ref{thm:A}, and Section~\ref{sec:counting} proves Theorems~\ref{thm:B} and~\ref{thm:C}.  Sections~\ref{sec:index4}--\ref{sec:trinomial} concern index four: Section~\ref{sec:index4} decomposes the count according to the vanishing order and proves Theorem~\ref{thm:D}~\ref{D:strata}, and Sections~\ref{sec:binomial} and~\ref{sec:trinomial} evaluate the sums $B(r)+B(r+s)$ and the trinomial counts and complete the proof of Theorem~\ref{thm:D}.

\section{Exact index and values at roots of unity}\label{sec:prelim}
We keep the notation of Section~\ref{subsec:main}.  Whenever $\operatorname{ind}(F)\mid\ell$, we use the presentation
\begin{equation}\label{eq:indexform}
 F(X)=\beta+X^ah(X^s),\qquad 1\leqslant a\leqslant s,\quad \deg h<\ell,
\end{equation}
of Section~\ref{subsec:index}.  The vanishing order of $F-\beta$ at $0$ is $r=a+ks$, where $k$ is the least exponent with a nonzero coefficient in $h$.

For $h(T)=\sum_{k=0}^{\ell-1}c_kT^k$ we call
\[
 S(h):=\{k\in\mathbb Z/\ell\mathbb Z:c_k\neq0\}
\]
the \emph{coefficient support} of $h$, and for a nonempty subset $S\subseteq\mathbb Z/\ell\mathbb Z$ we put
\[
 g_\ell(S):=\gcd\bigl(\ell,\{k-k':k,k'\in S\}\bigr),
\]
so that $g_\ell(S)=\ell$ when $\lvert S\rvert=1$.

\begin{proposition}[Index and coefficient support]\label{prop:support-gcd}
Let $F$ be as in \eqref{eq:indexform} with $h\neq0$.  Then
\[
 \operatorname{ind}(F)=\frac{\ell}{g_\ell(S(h))}.
\]
In particular, $F$ has exact index $\ell$ if and only if $g_\ell(S(h))=1$.
\end{proposition}

\begin{proof}
The nonzero exponents of $F-\beta$ are the numbers $a+ks$ with $k\in S(h)$.  Their pairwise differences are $s(k-k')$, and $q-1=\ell s$.  Hence $\Delta(F)=s\,g_\ell(S(h))$, which gives the formula.
\end{proof}

In particular, if $h$ has a single nonzero coefficient, then $F$ has index $1$.  Thus for $\ell\geqslant2$ no monomial has exact index $\ell$, while for $\ell=1$ every polynomial of index $1$ has a single nonzero coefficient in $h$.

For $h$ as above we write
\[
 A_j:=h(\zeta^j)\qquad(j\in\mathbb Z/\ell\mathbb Z)
\]
for its values on $\mu_\ell$.  Since $\ell\mid q-1$, the characteristic of $\mathbb F_q$ does not divide $\ell$, and the coefficients are recovered by the inverse discrete Fourier transform
\begin{equation}\label{eq:IDFT}
 c_k=\frac1\ell\sum_{j=0}^{\ell-1}A_j\zeta^{-jk}.
\end{equation}

\begin{lemma}[Cyclotomic criterion]\label{lem:cyclo}
Let $e\geqslant1$ and $h\in\mathbb F_q[T]$.  Then $X^eh(X^s)$ permutes $\mathbb F_q$ if and only if $\gcd(e,s)=1$ and the map $u\mapsto u^eh(u)^s$ permutes $\mu_\ell$.  In particular, $F$ as in \eqref{eq:indexform} is a PP if and only if $\gcd(a,s)=1$, all $A_j$ are nonzero, and the elements $\zeta^{aj}A_j^s$ $(0\leqslant j<\ell)$ are the $\ell$ distinct elements of $\mu_\ell$.
\end{lemma}

\begin{proof}
For $1\leqslant e\leqslant s$ and $\deg h<\ell$ this is the standard cyclotomic reduction \cite[Lemma~2.1 and Theorem~2.2]{AkbaryWang2007}: injectivity on the fibers of $x\mapsto x^s$ is equivalent to $\gcd(e,s)=1$, and the map induced on $\mu_\ell$ is $u\mapsto u^eh(u)^s$.  In general, write $e=a'+ks$ with $1\leqslant a'\leqslant s$ and let $h_1$ be the remainder of $T^kh(T)$ upon division by $T^\ell-1$.  The polynomials $X^eh(X^s)$ and $X^{a'}h_1(X^s)$ both vanish at $0$ and agree at every $x\neq0$, because $x^s\in\mu_\ell$; moreover $\gcd(e,s)=\gcd(a',s)$ and $u^eh(u)^s=u^{a'}h_1(u)^s$ for $u\in\mu_\ell$.  This reduces the general case to the previous one.  For the last assertion, note that adding the constant $\beta$ does not affect bijectivity, and that for $u\in\mu_\ell$ the element $u^ah(u)^s$ is $0$ if $h(u)=0$ and lies in $\mu_\ell$ otherwise, since $h(u)^{s\ell}=h(u)^{q-1}$.
\end{proof}

For a PP $F$ as in \eqref{eq:indexform}, the induced permutation $\sigma=\sigma_F$ of Section~\ref{subsec:main} is therefore characterized by
\[
 \zeta^{aj}A_j^s=\zeta^{\sigma(j)}\qquad(j\in\mathbb Z/\ell\mathbb Z).
\]
Prescribing the values $A_j^s\in\mu_\ell$ and applying \eqref{eq:IDFT} is the classical way of generating such polynomials.  The next theorem identifies, among the vectors $(A_j)$, those that give exact index $\ell$.

\begin{theorem}[Fourier criterion for exact index]\label{thm:fourier-exact}
Let $h\neq0$ have degree less than $\ell$.  The following are equivalent:
\begin{enumerate}[label=\textup{(\roman*)}]
\item\label{fe:index} $\operatorname{ind}(X^ah(X^s))=\ell$;
\item\label{fe:gcd} $g_\ell(S(h))=1$;
\item\label{fe:eigen} for every prime $p\mid\ell$ and every $b\in\mathbb Z/p\mathbb Z$, the identities
\begin{equation}\label{eq:fourier-obstruction}
 A_{j+\ell/p}=\zeta^{b\ell/p}A_j\qquad(j\in\mathbb Z/\ell\mathbb Z)
\end{equation}
do not all hold.
\end{enumerate}
\end{theorem}

\begin{proof}
The equivalence of \ref{fe:index} and \ref{fe:gcd} is Proposition~\ref{prop:support-gcd}, and \ref{fe:gcd} holds if and only if no prime $p\mid\ell$ divides $g_\ell(S(h))$.  If $p\mid g_\ell(S(h))$, then all $k\in S(h)$ lie in one residue class $b$ modulo $p$, and therefore
\[
 A_{j+\ell/p}=\sum_{k\in S(h)}c_k\zeta^{jk}\zeta^{(\ell/p)k}=\zeta^{b\ell/p}A_j.
\]
Conversely, if \eqref{eq:fourier-obstruction} holds, then by \eqref{eq:IDFT} each $c_k$ satisfies $c_k=\zeta^{(\ell/p)(b-k)}c_k$.  Since $\zeta^{\ell/p}$ has order $p$, a coefficient $c_k$ can be nonzero only when $k\equiv b\pmod p$.  Thus $p\mid g_\ell(S(h))$.
\end{proof}

For $\ell=1$ the condition in \ref{fe:eigen} is empty, in agreement with the remark after Proposition~\ref{prop:support-gcd}.  For $\ell=2$ the excluded vectors form the two eigenspaces $A_{j+1}=A_j$ and $A_{j+1}=-A_j$, which correspond under \eqref{eq:IDFT} to the two singleton supports.  For $\ell=3$ the excluded vectors are those of the form $A_j=\zeta^{bj}A_0$, which correspond to monomials.  For $\ell=4$ the excluded vectors form the union of the two eigenspaces $A_{j+2}=A_j$ and $A_{j+2}=-A_j$.

\section{The orbit structure}\label{sec:orbits}
In this section $F=\beta+X^ah(X^s)$ is a PP written as in \eqref{eq:indexform}, with values $A_j=h(\zeta^j)$ and induced permutation $\sigma$.

\subsection{Normalization and the effect of QM transformations}

\begin{lemma}[Normalization of the residue exponent]\label{lem:residue-one}
There is a positive integer $d$ with $\gcd(d,q-1)=1$ such that $F(X^d)$ has residue exponent $1$.  In particular, every QM class of PPs of exact index $\ell$ contains a polynomial with residue exponent $1$.
\end{lemma}

\begin{proof}
By Lemma~\ref{lem:cyclo}, $\gcd(a,s)=1$; choose $d_0$ with $d_0a\equiv1\pmod s$.  We lift $d_0$ to a unit modulo $\ell s$.  For every prime $p\mid\ell$ with $p\nmid s$, impose in addition $d\equiv1\pmod p$; together with $d\equiv d_0\pmod s$ these congruences are compatible by the Chinese remainder theorem.  Primes dividing both $\ell$ and $s$ do not divide $d$, because $d$ is a unit modulo $s$.  Thus $\gcd(d,\ell s)=1$ and $da\equiv1\pmod s$, so the substitution $X\mapsto X^d$ produces a polynomial with residue exponent $1$.  Since QM equivalence preserves the permutation property and the index, the last assertion follows.
\end{proof}

\begin{lemma}[Effect of a QM transformation]\label{lem:qm-action}
Let $u,b\in\mathbb F_q^*$ and let $d$ be a positive integer with $\gcd(d,q-1)=1$.  Write $G\equiv uF(bX^d)$ as $G=u\beta+X^{a'}h'(X^s)$ in the form \eqref{eq:indexform}, and put
\[
 \kappa:=\frac{da-a'}{s},\qquad \delta:=d\bmod\ell,
\]
and let $t\in\mathbb Z/\ell\mathbb Z$ be defined by $b^s=\zeta^t$.  Then
\[
 h'(T)\equiv ub^aT^\kappa h(\zeta^tT^\delta)\pmod{T^\ell-1}.
\]
Consequently, the values $A'_j=h'(\zeta^j)$ satisfy $A'_j=ub^a\zeta^{\kappa j}A_{t+\delta j}$, and the normalized values $U_j:=A_j/A_0$ and the induced permutation $\sigma'$ of $G$ transform by
\begin{equation}\label{eq:QM-action}
 U'_j=\zeta^{\kappa j}\frac{U_{t+\delta j}}{U_t},\qquad
 \sigma'(j)-\sigma'(0)=\sigma(t+\delta j)-\sigma(t).
\end{equation}
These formulas depend only on $\kappa\bmod\ell$, $\delta$ and $t$; in particular, choices of $b$ with the same $b^s$ induce the same transformation.
\end{lemma}

\begin{proof}
Since $a'\equiv da\pmod s$, $\kappa$ is a nonnegative integer.  For $x\in\mathbb F_q^*$ put $w:=x^s\in\mu_\ell$.  Then $x^{ds}=w^d=w^\delta$ and $x^{da}=x^{a'}w^\kappa$, so
\[
 uF(bx^d)-u\beta=ub^ax^{da}h(b^sx^{ds})=x^{a'}\cdot ub^aw^\kappa h(\zeta^tw^\delta).
\]
Both sides vanish at $x=0$.  A polynomial of degree less than $\ell$ is determined by its values on $\mu_\ell$, and $a'+s(\ell-1)\leqslant q-1$; hence the reduced representative of $G$ is $u\beta+X^{a'}h'(X^s)$ with $h'$ as stated, and $A'_j=ub^a\zeta^{\kappa j}A_{t+\delta j}$.  Dividing by $A'_0=ub^aA_t$ gives the first formula in \eqref{eq:QM-action}.  For the second, use $b^{as}=\zeta^{at}$, $A_{t+\delta j}^s=\zeta^{\sigma(t+\delta j)-a(t+\delta j)}$, $a'+\kappa s=da$ and $d\equiv\delta\pmod\ell$:
\[
 \zeta^{a'j}(A'_j)^s=u^s\zeta^{at}\zeta^{daj}\zeta^{\sigma(t+\delta j)-at-a\delta j}=u^s\zeta^{\sigma(t+\delta j)}.
\]
Since $u^s\in\mu_\ell$, the induced permutation of $G$ is $j\mapsto\sigma(t+\delta j)$ followed by a translation, which gives the second formula.
\end{proof}

\subsection{Zero constant term: a finite quotient}
Assume now that $F(0)=0$.  By Lemma~\ref{lem:residue-one} and output scaling, every QM class of PPs of exact index $\ell$ vanishing at $0$ contains a polynomial $Xh(X^s)$ with $h(1)=1$, that is, with $a=1$ and $A_0=1$.  We call such polynomials \emph{normalized}.  For a normalized PP one has $\sigma(0)=0$ and $U_j=A_j$, and
\[
 U_j^s=\zeta^{\sigma(j)-j}\qquad(j\in\mathbb Z/\ell\mathbb Z).
\]

If two normalized polynomials are QM equivalent, say $G\equiv uF(bX^d)$, then $a=a'=1$ forces $d\equiv1\pmod s$.  Let
\[
 \mathcal D_{\ell,s}:=\{d\in(\mathbb Z/\ell s\mathbb Z)^*:d\equiv1\pmod s\}.
\]
Every $d\in\mathcal D_{\ell,s}$ has a unique expression $d\equiv1+\kappa s\pmod{\ell s}$ with $\kappa\in\mathbb Z/\ell\mathbb Z$, which agrees with the $\kappa$ of Lemma~\ref{lem:qm-action}, and then $\delta\equiv1+\kappa s\pmod\ell$.  Conversely, for $\kappa\in\mathbb Z/\ell\mathbb Z$ the class $1+\kappa s$ is a unit modulo $\ell s$ exactly when $\delta$ is a unit modulo $\ell$, since it is $\equiv1\pmod s$; we call such $\kappa$ admissible.  Let $\Gamma^0_{\ell,s}$ be the set of pairs $\omega=(\kappa,t)$ with $\kappa$ admissible and $t\in\mathbb Z/\ell\mathbb Z$; the pair $(\kappa,t)$ records the substitutions $X\mapsto bX^d$ with $d\equiv1+\kappa s\pmod{\ell s}$ and $b^s=\zeta^t$.  Since composing $X\mapsto b_1X^{d_1}$ with $X\mapsto b_2X^{d_2}$ gives $X\mapsto b_1b_2^{d_1}X^{d_1d_2}$, these pairs form a group, isomorphic via $(\kappa,t)\mapsto(t,1+\kappa s)$ to the semidirect product $(\mathbb Z/\ell\mathbb Z)\rtimes\mathcal D_{\ell,s}$ with multiplication $(t_1,d_1)(t_2,d_2)=(t_1+d_1t_2,\,d_1d_2)$.  By Lemma~\ref{lem:qm-action}, $\Gamma^0_{\ell,s}$ acts on the right on vectors $U=(U_j)$ with $U_0=1$ by the first formula in \eqref{eq:QM-action}.

For $\sigma\in S_\ell$ with $\sigma(0)=0$, let $\mathcal X^{\rm ex}_{\ell,s}(\sigma)$ be the set of vectors $U=(1,U_1,\dots,U_{\ell-1})$ satisfying $U_j^s=\zeta^{\sigma(j)-j}$ for all $j$ and the exact-index condition of Theorem~\ref{thm:fourier-exact}, and put
\[
 \mathcal X^{\rm ex}_{\ell,s}:=\bigsqcup_{\sigma(0)=0}\mathcal X^{\rm ex}_{\ell,s}(\sigma).
\]

\begin{theorem}[Finite quotient and canonical forms]\label{thm:zero-quotient}
Sending a normalized PP to its pair $(\sigma;U)$ is a bijection from the normalized PPs of exact index $\ell$ onto $\mathcal X^{\rm ex}_{\ell,s}$, and it identifies the QM classes of PPs of exact index $\ell$ vanishing at $0$ with
\[
 \mathcal X^{\rm ex}_{\ell,s}/\Gamma^0_{\ell,s}.
\]
After fixing $\gamma$, every orbit has a unique lexicographically least representative in the coordinates $(\sigma(0),\dots,\sigma(\ell-1);\log_\gamma U_1,\dots,\log_\gamma U_{\ell-1})$, where $\log_\gamma U_j\in\{0,\dots,q-2\}$ denotes the discrete logarithm to the base $\gamma$.
\end{theorem}

\begin{proof}
Given $(\sigma;U)\in\mathcal X^{\rm ex}_{\ell,s}$, the inverse Fourier transform \eqref{eq:IDFT} gives the unique polynomial $h$ of degree less than $\ell$ with $h(\zeta^j)=U_j$, namely
\[
 h(T)=\sum_{k=0}^{\ell-1}c_kT^k,\qquad c_k=\frac1\ell\sum_{j=0}^{\ell-1}U_j\zeta^{-jk}.
\]
The defining equations give the distinct induced values required by Lemma~\ref{lem:cyclo}, so $Xh(X^s)$ is a PP; it is normalized because $h(1)=U_0=1$, and it has exact index $\ell$ by Theorem~\ref{thm:fourier-exact}.  Together with the discussion above, this proves the first assertion and shows that every class meets the set of normalized polynomials.

If two normalized PPs are QM equivalent, say $G\equiv uF(bX^d)$, then $d\equiv1\pmod s$.  Write $d\equiv1+\kappa s\pmod{\ell s}$ and $b^s=\zeta^t$.  Requiring the transformed value at $1$ to remain $1$ forces $u=(bA_t)^{-1}$, so by Lemma~\ref{lem:qm-action} the vectors are related by $\omega=(\kappa,t)$.  Conversely, every $\omega=(\kappa,t)\in\Gamma^0_{\ell,s}$ is realized: with $b=\gamma^t$, $d\equiv1+\kappa s\pmod{\ell s}$ and $u=(bA_t)^{-1}$, the polynomial $uF(bX^d)$ is normalized.  On coefficients this transformation reads
\[
 h'(T)\equiv A_t^{-1}T^\kappa h(\zeta^tT^\delta)\pmod{T^\ell-1}.
\]
Since $\delta$ is a unit modulo $\ell$, it sends the coefficient support to $\kappa+\delta S(h)$ and preserves $g_\ell$, in accordance with the invariance of the index.  Hence the QM classes correspond to the $\Gamma^0_{\ell,s}$-orbits on $\mathcal X^{\rm ex}_{\ell,s}$.  Each orbit is finite, and a point is determined by its coordinates; this gives the unique lexicographically least representative.
\end{proof}

\subsection{Double cosets and stabilizers}
The group $\mathcal A_{\ell,s}=(\mathbb Z/\ell\mathbb Z)\rtimes H_{\ell,s}$ of Section~\ref{subsec:main} acts on $\mathbb Z/\ell\mathbb Z$ by $j\mapsto t+\delta j$.

\begin{lemma}\label{lem:pi}
Reduction modulo $\ell$ maps $\mathcal D_{\ell,s}$ onto $H_{\ell,s}$, with kernel consisting of the $g$ classes $d\equiv1\pmod{\operatorname{lcm}(\ell,s)}$ modulo $\ell s$, and $\lvert H_{\ell,s}\rvert=\varphi(\ell)/\varphi(g)$.  Consequently
\[
 \pi\colon\Gamma^0_{\ell,s}\longrightarrow\mathcal A_{\ell,s},\qquad(\kappa,t)\longmapsto(t,1+\kappa s)
\]
is a surjective homomorphism with $\lvert\ker\pi\rvert=g$, and
\[
 \lvert\Gamma^0_{\ell,s}\rvert=\ell g\,\frac{\varphi(\ell)}{\varphi(g)}.
\]
\end{lemma}

\begin{proof}
The simultaneous congruences $d\equiv1\pmod s$ and $d\equiv\delta\pmod\ell$ are solvable exactly when $\delta\equiv1\pmod g$, and if $\delta$ is a unit modulo $\ell$, every solution is a unit modulo $\ell s$.  The solutions with $\delta=1$ are the $g$ classes $d\equiv1\pmod{\operatorname{lcm}(\ell,s)}$ modulo $\ell s$.  Reduction $(\mathbb Z/\ell\mathbb Z)^*\to(\mathbb Z/g\mathbb Z)^*$ is onto, so $\lvert H_{\ell,s}\rvert=\varphi(\ell)/\varphi(g)$.  The multiplication rule of $\Gamma^0_{\ell,s}$ reduces modulo $\ell$ to the composition of affine maps, so $\pi$ is a homomorphism; its kernel consists of the pairs with $t=0$ and $1+\kappa s\equiv1\pmod\ell$, of which there are $g$.  Finally $\lvert\Gamma^0_{\ell,s}\rvert=\lvert\mathcal A_{\ell,s}\rvert\cdot g=\ell\lvert H_{\ell,s}\rvert g$.
\end{proof}

\begin{theorem}[Double cosets and stabilizers]\label{thm:double-coset}
Under the second formula in \eqref{eq:QM-action}, the $\Gamma^0_{\ell,s}$-orbits on the set of permutations $\sigma\in S_\ell$ with $\sigma(0)=0$ are indexed by
\[
 \mathcal T_\ell\backslash S_\ell/\mathcal A_{\ell,s};
\]
for some double cosets the corresponding sets $\mathcal X^{\rm ex}_{\ell,s}(\sigma)$ may be empty.  For $\sigma\in S_\ell$ with $\sigma(0)=0$, the stabilizer of $\sigma$ in $\Gamma^0_{\ell,s}$ is
\[
 K_\sigma=\pi^{-1}\bigl(\mathcal A_{\ell,s}\cap\sigma^{-1}\mathcal T_\ell\sigma\bigr),
\]
and consequently
\[
 \mathcal X^{\rm ex}_{\ell,s}/\Gamma^0_{\ell,s}
 \cong
 \bigsqcup_{[\sigma]\in\mathcal T_\ell\backslash S_\ell/\mathcal A_{\ell,s}}
 \mathcal X^{\rm ex}_{\ell,s}(\sigma)/K_\sigma.
\]
If all cycles of $\alpha\in\mathcal A_{\ell,s}$ have a common length $L$, say $m=\ell/L$ cycles, then $\alpha$ fixes exactly
\[
 \varphi(L)L^{m-1}(m-1)!
\]
induced permutations $\sigma$ with $\sigma(0)=0$; if its cycle lengths are not all equal, it fixes none.
\end{theorem}

\begin{proof}
Let $\Omega:=\{\sigma\in S_\ell:\sigma(0)=0\}$.  The map $\Omega\to\mathcal T_\ell\backslash S_\ell$ sending $\sigma$ to its left coset $\mathcal T_\ell\sigma$ is a bijection, because every left coset contains exactly one permutation fixing $0$.  Under this identification, the second formula in \eqref{eq:QM-action} is the right action of $\pi(\omega)\in\mathcal A_{\ell,s}$ on $\mathcal T_\ell\backslash S_\ell$.  Its orbit set is therefore $\mathcal T_\ell\backslash S_\ell/\mathcal A_{\ell,s}$.

For $\sigma\in\Omega$, an element $\alpha\in\mathcal A_{\ell,s}$ fixes $\sigma$ exactly when $\mathcal T_\ell\sigma\alpha=\mathcal T_\ell\sigma$, that is, when $\alpha\in\mathcal A_{\ell,s}\cap\sigma^{-1}\mathcal T_\ell\sigma$.  Pulling this stabilizer back through the surjection $\pi$ gives $K_\sigma$.  The projection $\mathcal X^{\rm ex}_{\ell,s}\to\Omega$ is equivariant.  Hence every $\Gamma^0_{\ell,s}$-orbit lying over a given double coset meets the fiber $\mathcal X^{\rm ex}_{\ell,s}(\sigma)$ over a chosen representative $\sigma$, and two points of this fiber lie in the same orbit exactly when they differ by an element of $K_\sigma$.  This gives the decomposition.

Now let $\alpha(j)=t+\delta j$.  A normalized $\sigma$ is fixed by $\alpha$ exactly when $\sigma(\alpha(j))-\sigma(t)=\sigma(j)$ for every $j$, that is, when $\sigma\alpha\sigma^{-1}$ is the translation by $\sigma(t)$.  Hence $\alpha$ must be conjugate to a translation.  A translation of additive order $L$ has cycle type $L^m$ with $m=\ell/L$, and there are $\varphi(L)$ such translations.  For a fixed target translation, the conjugators form a coset of its centralizer, a group of order $L^mm!$ that is transitive on the $\ell=Lm$ points; therefore exactly a fraction $1/\ell$ of these conjugators satisfy $\sigma(0)=0$.  Thus each target translation contributes $L^{m-1}(m-1)!$ normalized conjugators, which gives the asserted count.
\end{proof}

\subsection{Nonzero constant term}\label{subsec:nonzero}
QM equivalence involves no output translations, so the PPs with $F(0)\neq0$ form a union of QM classes.  If $F(0)=\beta\neq0$, multiplication by $\beta^{-1}$ normalizes the constant term to $1$, say $\beta^{-1}F=1+f$ with $f(0)=0$.  Output translations and scalings preserve bijectivity and the nonconstant exponent support.  Hence $f\mapsto1+f$ is a bijection from the PPs of exact index $\ell$ vanishing at $0$ to the PPs of exact index $\ell$ with constant term $1$.  If two polynomials with constant term $1$ are QM equivalent, evaluation at $0$ forces the output scalar to be $1$, so they differ by a substitution $X\mapsto vX^d$.

\begin{lemma}[Freeness]\label{lem:free-input}
The group
\[
 \mathcal M_q:=\{X\mapsto vX^d:v\in\mathbb F_q^*,\ \gcd(d,q-1)=1\}
\]
of input substitutions acts freely on the PPs with constant term $1$, and $\lvert\mathcal M_q\rvert=(q-1)\varphi(q-1)$.
\end{lemma}

\begin{proof}
Write a PP with constant term $1$ as $F=1+f$; then $f$ is a PP vanishing at $0$.  If $F(vX^d)=F(X)$, then $f(vX^d)=f(X)$.  Since $f$ is bijective, $vx^d=x$ for all $x\in\mathbb F_q$; evaluation at $1$ and at a generator of $\mathbb F_q^*$ gives $v=1$ and $d\equiv1\pmod{q-1}$.
\end{proof}

\begin{theorem}[Classification for nonzero constant term]\label{thm:nonzero-structure}
Let $F=1+Xh(X^s)$ be a PP of exact index $\ell$ with constant term $1$ and residue exponent $1$, normalized by $A_0=\gamma^{\sigma(0)}$, and put $U_j:=A_j/A_0$; here and below, in an exponent of $\gamma$ an element of $\mathbb Z/\ell\mathbb Z$ is represented by its integer representative in $\{0,\dots,\ell-1\}$.  Every $\mathcal M_q$-orbit of PPs of exact index $\ell$ with constant term $1$ contains such polynomials, and for them
\[
 U_j^s=\zeta^{\sigma(j)-j-\sigma(0)}.
\]
An element $\omega=(\kappa,t)$ of $\Gamma^0_{\ell,s}$, with $\delta\equiv1+\kappa s\pmod\ell$, acts on these data by
\[
 \sigma'(j)=\sigma(t+\delta j),\qquad
 U'_j=\zeta^{\kappa j}\frac{U_{t+\delta j}}{U_t},
\]
so the action on induced permutations is the right action of $\mathcal A_{\ell,s}$ on $S_\ell$.  Hence every QM class of PPs of exact index $\ell$ with constant term $1$ projects to a right orbit in
\[
 S_\ell/\mathcal A_{\ell,s},
\]
and only those right orbits occur whose coefficient fibers contain a point of exact index $\ell$.  The right action of $\mathcal A_{\ell,s}$ on $S_\ell$ is free; hence, for every induced permutation that occurs, the stabilizer in $\Gamma^0_{\ell,s}$ is $\ker\pi$, of order $g$, and $\ker\pi$ acts freely on each coefficient fiber.  After fixing $\gamma$ and taking discrete logarithms $\log_\gamma U_j\in\{0,\dots,q-2\}$, every class has a unique lexicographically least normalized representative in the coordinates
\[
 (\sigma(0),\dots,\sigma(\ell-1);\log_\gamma U_1,\dots,\log_\gamma U_{\ell-1}).
\]
\end{theorem}

\begin{proof}
By Lemma~\ref{lem:residue-one} we may assume $a=1$.  Since $A_0^s=\zeta^{\sigma(0)}$, the scalar $b_0:=\gamma^{\sigma(0)}/A_0$ satisfies $b_0^s=1$; the substitution $X\mapsto b_0X$ preserves $a=1$ and gives $A_0=\gamma^{\sigma(0)}$.  Thus every orbit of input substitutions meets the set of normalized polynomials, and the formula for $U_j^s$ follows from $\zeta^jA_j^s=\zeta^{\sigma(j)}$.

If two normalized polynomials differ by $X\mapsto bX^d$, then $d\equiv1\pmod s$.  Write $d\equiv1+\kappa s\pmod{\ell s}$, $\delta\equiv d\pmod\ell$ and $b^s=\zeta^t$.  For the source $Xh(X^s)$, substitution gives $h'(T)\equiv bT^\kappa h(\zeta^tT^\delta)\pmod{T^\ell-1}$, hence $A'_j=b\zeta^{\kappa j}A_{t+\delta j}$.  The normalization of the target forces $b=\gamma^{\sigma(t)}/A_t$, whose $s$-th power is indeed $\zeta^t$; taking $s$-th powers, and dividing by $A'_0$, gives the displayed action.  Conversely this choice of $b$ realizes every $(\kappa,t)$, so the orbits of input substitutions on normalized polynomials are exactly the $\Gamma^0_{\ell,s}$-orbits.

The projection to $S_\ell$ is equivariant for the right action of $\mathcal A_{\ell,s}$, which is free; hence stabilizers lie in $\ker\pi$.  A nonidentity element of $\ker\pi$ has $t=0$, $\delta=1$ and $\kappa\neq0$, and it sends $U_1$ to $\zeta^\kappa U_1\neq U_1$.  Thus $\ker\pi$ acts freely on each fiber; the right orbits that occur need not exhaust $S_\ell/\mathcal A_{\ell,s}$.  The coordinates determine the point, which gives the unique lexicographically least representative in each finite orbit.
\end{proof}

\begin{proof}[Proof of Theorem~\textup{\ref{thm:A}}]
Part \ref{A:normal} follows from Lemma~\ref{lem:residue-one}, Theorem~\ref{thm:zero-quotient} and Lemma~\ref{lem:pi}.  Part \ref{A:double} follows from Lemma~\ref{lem:pi}, the second formula in \eqref{eq:QM-action} and Theorem~\ref{thm:double-coset}, and part \ref{A:fixed} is the last assertion of Theorem~\ref{thm:double-coset}, since a normalized $\sigma$ corresponds to the left coset $\mathcal T_\ell\sigma$.  Part \ref{A:nonzero} is the discussion at the beginning of Section~\ref{subsec:nonzero} together with Lemma~\ref{lem:free-input}.  The canonical representatives are provided by Theorems~\ref{thm:zero-quotient} and~\ref{thm:nonzero-structure}.
\end{proof}

\section{Counting the classes}\label{sec:counting}
Let $\mu$ and $\varphi$ be the M\"obius and Euler functions, and let $\mathfrak Q_\ell(q)$, $\mathfrak R_\ell(q)$ and $\mathfrak A_\ell(q)$ be as in Section~\ref{subsec:main}.

\subsection{Fixed points and cycle congruences}
Fix $\sigma\in S_\ell$ with $\sigma(0)=0$ and $\omega=(\kappa,t)\in K_\sigma$, and let $\alpha:=\pi(\omega)$, so that $\alpha(j)=t+\delta j$ with $\delta\equiv1+\kappa s\pmod\ell$.  Choose $\bar e_j\in\{0,\dots,\ell-1\}$ with $\bar e_j\equiv\sigma(j)-j\pmod\ell$, and let $\mathcal X_{\ell,s}(\sigma)$ be the set of all vectors $U$ with $U_0=1$ and $U_j^s=\zeta^{\sigma(j)-j}$ for all $j$.  Every $U\in\mathcal X_{\ell,s}(\sigma)$ can be written uniquely as $U_j=\gamma^{\bar e_j+\ell x_j}$ with $x_j\in\mathbb Z/s\mathbb Z$ and $x_0=0$.  Let $d\mid\ell$, put $m:=\ell/d$, and let $c\in\mathbb Z/d\mathbb Z$.  The vector $U$ is fixed by $\omega$ exactly when
\begin{equation}\label{eq:fixed-linear}
 x_j-x_{\alpha(j)}+x_t\equiv\frac{\kappa sj+\bar e_{\alpha(j)}-\bar e_t-\bar e_j}{\ell}\pmod s\qquad(j\in\mathbb Z/\ell\mathbb Z),
\end{equation}
where the numerator is divisible by $\ell$ because $\omega\in K_\sigma$.  As in the proof of Theorem~\ref{thm:fourier-exact}, the coefficient support of the corresponding $h$ lies in the class $c$ modulo $d$ exactly when $U_{j+m}=\zeta^{cm}U_j$ for all $j$, that is, when
\begin{equation}\label{eq:support-linear}
 x_{j+m}-x_j\equiv\frac{scm-(\bar e_{j+m}-\bar e_j)}{\ell}\pmod s\qquad(j\in\mathbb Z/\ell\mathbb Z);
\end{equation}
this system has no solution unless every numerator is divisible by $\ell$.

\begin{lemma}[Cycle count for a fixed-point term]\label{lem:cycle-count}
Let $\nu_{\omega,\sigma;d,c}$ be the number of solutions $(x_j)$ of \eqref{eq:fixed-linear} and \eqref{eq:support-linear} with $x_0=0$.  Put $\bar\alpha(j):=t+\delta j\pmod m$ and $\bar t:=t\pmod m$.  Then $\nu_{\omega,\sigma;d,c}=0$ unless
\begin{equation}\label{eq:cycle-compat}
 \bar e_{j+m}-\bar e_j\equiv scm\pmod\ell\quad(j\in\mathbb Z/\ell\mathbb Z),\qquad d\mid\kappa(1+cs).
\end{equation}
If these conditions hold, choose $\tilde\kappa\in\{0,\dots,\ell-1\}$ with $\tilde\kappa\equiv\kappa(1+cs)\pmod\ell$, and integers $0\leqslant\rho_j<\ell$ with $\rho_j\equiv\bar e_j-csj\pmod\ell$ for $0\leqslant j<m$, and set, identifying $\mathbb Z/m\mathbb Z$ with $\{0,\dots,m-1\}$,
\[
 b_j:=\frac{s\tilde\kappa j-\rho_j+\rho_{\bar\alpha(j)}-\rho_{\bar t}}{\ell}\in\mathbb Z\qquad(j\in\mathbb Z/m\mathbb Z).
\]
Let $C_1,\dots,C_N$ be the cycles of $\bar\alpha$ on $\mathbb Z/m\mathbb Z$, of lengths $L_1,\dots,L_N$, and put $\lambda:=\gcd(L_1,\dots,L_N)$ and $B_i:=\sum_{j\in C_i}b_j$.  Then
\[
 \nu_{\omega,\sigma;d,c}=\begin{cases}
 s^{N-1}\gcd(\lambda,s),& \text{if the congruences }L_iy\equiv B_i\pmod s\ (1\leqslant i\leqslant N)\\
 &\text{have a common solution }y\in\mathbb Z/s\mathbb Z,\\
 0,&\text{otherwise.}
 \end{cases}
\]
Moreover, the homogeneous integer system associated with \eqref{eq:fixed-linear} and \eqref{eq:support-linear} presents the group $\mathbb Z^{N-1}\oplus\mathbb Z/\lambda\mathbb Z$.
\end{lemma}

\begin{proof}
Write $V_j:=\gamma^{\bar e_j+\ell x_j}$ and $V'_j:=\zeta^{-cj}V_j$.  The support equations say that $V'_{j+m}=V'_j$, while the fixed-point equations become
\[
 V'_j=\zeta^{\kappa(1+cs)j}\frac{V'_{t+\delta j}}{V'_t},
\]
because $\kappa+c(\delta-1)\equiv\kappa(1+cs)\pmod\ell$.  Replacing $j$ by $j+m$ and using periodicity forces $\zeta^{\kappa(1+cs)m}=1$, that is, $d\mid\kappa(1+cs)$; under this condition the equations reduce to $V'_j=\zeta^{\tilde\kappa j}V'_{\bar\alpha(j)}/V'_{\bar t}$ on $\mathbb Z/m\mathbb Z$.  Periodicity of $(V'_j)^s=\zeta^{\bar e_j-csj}$ is the first condition in \eqref{eq:cycle-compat}.  Writing $V'_j=\gamma^{\rho_j+\ell z_j}$ with $z_0=0$, we obtain
\[
 z_j-z_{\bar\alpha(j)}+z_{\bar t}\equiv b_j\pmod s,
\]
and the solutions $(x_j)$ correspond bijectively to the solutions $(z_j)_{j\in\mathbb Z/m\mathbb Z}$.

To see that $b_j\in\mathbb Z$, put $k_j:=(t+\delta j-\bar\alpha(j))/m$ and $k_0:=(t-\bar t)/m$, so that $m(k_j-k_0)=\delta j-\bar\alpha(j)+\bar t$.  Iterating the first congruence in \eqref{eq:cycle-compat} gives
\[
 \bar e_{t+\delta j}-\bar e_t\equiv\bar e_{\bar\alpha(j)}-\bar e_{\bar t}+(k_j-k_0)scm\pmod\ell.
\]
Combining this with the stabilizer relation $\kappa sj+\bar e_{t+\delta j}-\bar e_t-\bar e_j\equiv0\pmod\ell$ and with $\tilde\kappa\equiv\kappa(1+cs)$ yields
\[
 s\tilde\kappa j-\rho_j+\rho_{\bar\alpha(j)}-\rho_{\bar t}\equiv csj(1+\kappa s-\delta)\equiv0\pmod\ell,
\]
since $\delta\equiv1+\kappa s\pmod\ell$.

Now $\bar\alpha(0)=\bar t$, so $0$ and $\bar t$ lie in the same cycle.  Setting $y:=z_{\bar t}$, the equations read $z_{\bar\alpha(j)}=z_j+y-b_j$, and summing around $C_i$ gives $L_iy\equiv B_i\pmod s$.  Conversely, each common solution $y$ determines the variables on the cycle containing $0$ from $z_0=0$, and leaves one free initial value on each of the other $N-1$ cycles.  There are therefore $s^{N-1}$ solutions for each admissible $y$, and the set of admissible $y$ is either empty or a coset of a subgroup of order $\gcd(s,L_1,\dots,L_N)=\gcd(\lambda,s)$.

For the homogeneous integer system, the support rows first identify variables whose indices are congruent modulo $m$.  The remaining relations are $z_{\bar\alpha(j)}=z_j+y$ with $y=z_{\bar t}$.  Since $\bar\alpha(0)=\bar t$, on the cycle through $0$ the value $z_0=0$ is fixed; every other cycle contributes one free generator, while closing the cycle $C_i$ imposes $L_iy=0$.  The ideal generated by the $L_i$ is $\lambda\mathbb Z$, so the presented group is $\mathbb Z^{N-1}\oplus\mathbb Z/\lambda\mathbb Z$.
\end{proof}

\subsection{The class number formula}

\begin{theorem}[Number of classes for arbitrary index]\label{thm:class-number}
For $\omega\in K_\sigma$ put
\[
 F^{\rm ex}_{\omega,\sigma}:=\sum_{d\mid\ell}\mu(d)\sum_{c\bmod d}\nu_{\omega,\sigma;d,c},
\]
with $\nu_{\omega,\sigma;d,c}$ as in Lemma~\textup{\ref{lem:cycle-count}}.  Then, choosing in each double coset a representative $\sigma$ with $\sigma(0)=0$,
\[
 \mathfrak Q_\ell(q)=\sum_{[\sigma]\in\mathcal T_\ell\backslash S_\ell/\mathcal A_{\ell,s}}\frac1{\lvert K_\sigma\rvert}\sum_{\omega\in K_\sigma}F^{\rm ex}_{\omega,\sigma}.
\]
Each summand is independent of the choice of $\sigma$, being the number of classes lying over the double coset by Theorem~\textup{\ref{thm:double-coset}}.  Thus the classification of the PPs of exact index $\ell$ vanishing at $0$, and the number of their QM classes, are finite and explicit for every $\ell\mid q-1$.  The contribution of the identity simplifies to
\begin{equation}\label{eq:exact-cardinality}
 \lvert\mathcal X^{\rm ex}_{\ell,s}\rvert=
 \sum_{d\mid\ell}\mu(d)\,
 \gcd(d,s)\frac{\varphi(d)}{\varphi(\gcd(d,s))}
 \left(\frac{\ell}{d}-1\right)!\,(ds)^{\ell/d-1}.
\end{equation}
\end{theorem}

\begin{proof}
For a nonempty support $S$, the condition $d\mid g_\ell(S)$ means that $S$ lies in a unique class $c\bmod d$, and by Fourier inversion this is exactly \eqref{eq:support-linear}.  Since $U_0=1$, the support of every point of $\mathcal X_{\ell,s}(\sigma)$ is nonempty.  Hence M\"obius inversion, in the form $\mathbf 1_{g_\ell(S)=1}=\sum_{d\mid g_\ell(S)}\mu(d)$, shows that $F^{\rm ex}_{\omega,\sigma}$ is the number of points of $\mathcal X^{\rm ex}_{\ell,s}(\sigma)$ fixed by $\omega$.  Burnside's lemma, applied to each quotient $\mathcal X^{\rm ex}_{\ell,s}(\sigma)/K_\sigma$ in Theorem~\ref{thm:double-coset}, gives the class number.

For \eqref{eq:exact-cardinality}, a normalized $h$ whose support lies in the class $c$ modulo $d$ can be written as $h(T)=T^c\tilde h(T^d)$ with $\deg\tilde h<\ell/d$.  By Lemma~\ref{lem:cyclo}, applied with $\ell/d$ in place of $\ell$, the polynomial $X^{1+cs}\tilde h(X^{ds})$ is a PP if and only if $\gcd(1+cs,d)=1$ and the induced map on $\mu_{\ell/d}$ is a permutation, since $\gcd(1+cs,s)=1$.  If $g'=\gcd(d,s)$, the Chinese remainder theorem gives $g'\varphi(d)/\varphi(g')$ admissible classes $c\bmod d$: at primes dividing $s$ all classes are allowed, and at primes $p\nmid s$ exactly $\varphi(p^a)$ classes modulo $p^a$ are allowed.  For $n=\ell/d$, each admissible class contributes $(n-1)!(ds)^{n-1}$ normalized points, from $(n-1)!$ normalized induced permutations and $ds$ choices for each remaining value.
\end{proof}

\subsection{Nonzero constant term}
Let $N_0(\ell,q)$ be the number of PPs of exact index $\ell$ vanishing at $0$.  Li, Guan and Wu \cite[Theorem~1.4 and its proof, pp.~841--842]{LiGuanWu2025} showed that
\begin{equation}\label{eq:LGW-general}
 N_0(\ell,q)=\sum_{d\mid\ell}\mu(\ell/d)\,d!\left(\frac{q-1}{d}\right)^d\varphi\left(\frac{q-1}{d}\right).
\end{equation}

\begin{proposition}[Number of classes with nonzero constant term]\label{prop:Rell}
\[
 \mathfrak R_\ell(q)=\frac{N_0(\ell,q)}{(q-1)\varphi(q-1)}.
\]
\end{proposition}

\begin{proof}
By Section~\ref{subsec:nonzero}, the PPs of exact index $\ell$ with nonzero constant term are the polynomials $\beta+f$ with $\beta\neq0$ and $f$ a PP of exact index $\ell$ vanishing at $0$; there are $(q-1)N_0(\ell,q)$ of them.  The group of QM transformations $(u,v,d)$ has $(q-1)^2\varphi(q-1)$ elements and acts freely on them: if $uF(vX^d)=F$ with $F(0)\neq0$, then $u=1$ by evaluation at $0$, and bijectivity gives $vX^d=X$, hence $v=1$ and $d\equiv1\pmod{q-1}$.  Dividing by the orbit size gives the formula.
\end{proof}

\subsection{Quasi-polynomiality}

\begin{lemma}[Cycle sums of affine permutations]\label{lem:affine-cycle-sum}
Let $m\geqslant1$ and let $\alpha$ be an affine permutation $j\mapsto t+\delta j$ of $\mathbb Z/m\mathbb Z$ with $\delta\in(\mathbb Z/m\mathbb Z)^*$.  Identify $\mathbb Z/m\mathbb Z$ with $\{0,\dots,m-1\}$.  Then every cycle $C$ of $\alpha$, of length $L$, satisfies
\begin{equation}\label{eq:affine-cycle-sum-divisibility}
 \gcd(m,L)\mid2\sum_{j\in C}j.
\end{equation}
\end{lemma}

\begin{proof}
We argue by induction on $m$, the case $m=1$ being trivial.  Let $m':=\gcd(m,L)$.  If $m'=m$, then $m\mid L$, so $L=m$ and $C$ is the full residue set; then $2\sum_{j\in C}j=m(m-1)$.  Otherwise $m'<m$, and reduction modulo $m'$ maps $\alpha$ to an affine permutation of $\mathbb Z/m'\mathbb Z$ and $C$ onto a cycle $C'$ of length $L'\mid L$, each point of $C'$ being the image of exactly $L/L'$ points of $C$.  Hence $2\sum_{j\in C}j\equiv(L/L')\,2\sum_{j\in C'}j\pmod{m'}$, and by induction $\gcd(m',L')\mid2\sum_{j\in C'}j$.  It remains to show that $m'\mid(L/L')\gcd(m',L')$.  For a prime $p$ write $x=v_p(m')$, $y=v_p(L)$ and $z=v_p(L')$; since $x,z\leqslant y$, one has $y-z+\min(x,z)\geqslant x$.  This gives \eqref{eq:affine-cycle-sum-divisibility}.
\end{proof}

\begin{theorem}[Quasi-polynomial counts]\label{thm:quasipolynomial}
Fix $\ell\geqslant2$ and let $P_\ell:=\ell$ for odd $\ell$ and $P_\ell:=2\ell$ for even $\ell$.  There are effectively computable polynomials $Q_{\ell,\varrho},R_{\ell,\varrho}\in\mathbb Q[X]$ $(\varrho\in\mathbb Z/P_\ell\mathbb Z)$ of degree at most $\ell-1$ such that, whenever $q=\ell s+1$ is a prime power and $s\equiv\varrho\pmod{P_\ell}$,
\[
 \mathfrak Q_\ell(q)=Q_{\ell,\varrho}(s),\qquad
 \mathfrak R_\ell(q)=R_{\ell,\varrho}(s),\qquad
 \mathfrak A_\ell(q)=Q_{\ell,\varrho}(s)+R_{\ell,\varrho}(s).
\]
Thus the periods of $\mathfrak Q_\ell$ and $\mathfrak A_\ell$ divide $P_\ell$, while the period of $\mathfrak R_\ell$ divides $\operatorname{rad}(\ell)$.  With $g=\gcd(\ell,s)$ as in Section~\ref{subsec:main}, the coefficient of $X^{\ell-1}$ in $Q_{\ell,\varrho}$ is
\[
 \frac{(\ell-1)!\,\varphi(g)}{\ell g\,\varphi(\ell)}.
\]
\end{theorem}

\begin{proof}
Fix a term $(\sigma,\omega,d,c)$ in Theorem~\ref{thm:class-number} and use the notation of Lemma~\ref{lem:cycle-count}.  On a fixed residue class of $s$ modulo $\ell$, the compatibility conditions \eqref{eq:cycle-compat}, the affine permutation $\bar\alpha$, its cycles $C_i$, and the data $\tilde\kappa$ and $\rho_j$ are constant; so are the groups $\mathcal A_{\ell,s}$, $\Gamma^0_{\ell,s}$ and $K_\sigma$.  Put $\Sigma_i:=\sum_{j\in C_i}j$, with $j\in\{0,\dots,m-1\}$.  Summing the definition of $b_j$ around $C_i$ gives
\begin{equation}\label{eq:cycle-sum-B}
 B_i(s)=\frac{s\tilde\kappa\Sigma_i-L_i\rho_{\bar t}}{\ell}.
\end{equation}
Choose integers $\theta_i$ with $\sum_i\theta_iL_i=\lambda$, and put $\Sigma:=\sum_i\theta_i\Sigma_i$ and $B^\ast(s):=\sum_i\theta_iB_i(s)$.  The congruences $L_iy\equiv B_i(s)\pmod s$ have a common solution exactly when
\begin{align}
 B_i(s)-\frac{L_i}{\lambda}B^\ast(s)&\equiv0\pmod s\qquad(1\leqslant i\leqslant N),\label{eq:cycle-free-compat}\\
 \gcd(\lambda,s)&\mid B^\ast(s).\label{eq:cycle-torsion-compat}
\end{align}
Indeed, any common solution satisfies $\lambda y\equiv B^\ast(s)$, which gives both conditions; conversely, \eqref{eq:cycle-torsion-compat} provides $y$ with $\lambda y\equiv B^\ast(s)$, and then \eqref{eq:cycle-free-compat} gives every cycle congruence.  By \eqref{eq:cycle-sum-B}, the left side of \eqref{eq:cycle-free-compat} equals
\[
 \frac{s\tilde\kappa}{\ell}\left(\Sigma_i-\frac{L_i}{\lambda}\Sigma\right),
\]
so its divisibility by $s$ is equivalent to the integrality of the constant $\frac{\tilde\kappa}{\ell}(\Sigma_i-\frac{L_i}{\lambda}\Sigma)$.  Thus the conditions \eqref{eq:cycle-free-compat} are constant on the chosen residue class.

Since $\lambda\mid\sum_iL_i=m\mid\ell$, the number $\gcd(\lambda,s)$ is unchanged by $s\mapsto s+\ell$, while \eqref{eq:cycle-sum-B} shows that $B^\ast(s)$ changes by $\tilde\kappa\Sigma$.  Lemma~\ref{lem:affine-cycle-sum}, applied to the cycles of $\bar\alpha$, gives $\lambda\mid\gcd(m,L_i)\mid2\Sigma_i$, hence $\lambda\mid2\Sigma$.  Thus two shifts by $\ell$ preserve \eqref{eq:cycle-torsion-compat}; if $\ell$ is odd, then $\lambda$ is odd and one shift suffices.  Whether the congruences $L_iy\equiv B_i(s)\pmod s$ have a common solution therefore depends only on $s\bmod P_\ell$.

By Lemma~\ref{lem:cycle-count}, each fixed-point term is $s^{N-1}\gcd(\lambda,s)$ times the indicator of a condition on $s\bmod P_\ell$.  Since $N\leqslant m\leqslant\ell$, the finite sum of Theorem~\ref{thm:class-number} is a quasi-polynomial of degree at most $\ell-1$.  A term of degree $\ell-1$ requires $N=m=\ell$, so $d=1$ and $\bar\alpha$ is the identity; hence $t=0$ and $\delta=1$.  The fixed-point equation then reads $U_j=\zeta^{\kappa j}U_j$, and $j=1$ gives $\kappa=0$.  Thus only the identity contributes in top degree.  Its contribution has leading term $(\ell-1)!\,s^{\ell-1}$, and by Lemma~\ref{lem:pi} $\lvert\Gamma^0_{\ell,s}\rvert=\ell g\varphi(\ell)/\varphi(g)$, which gives the stated coefficient.

Finally, by Proposition~\ref{prop:Rell} and \eqref{eq:LGW-general},
\[
 \mathfrak R_\ell(q)=\sum_{d\mid\ell}\mu(\ell/d)\,\frac{d!\,(\ell/d)^d s^{d-1}}{\ell}\cdot\frac{\varphi(\ell s/d)}{\varphi(\ell s)},
 \qquad
 \frac{\varphi(\ell s/d)}{\varphi(\ell s)}
 =\frac1d\prod_{\substack{p\mid d,\ v_p(d)=v_p(\ell)\\ p\nmid s}}\frac{p}{p-1},
\]
and the last factor depends on $s$ only through the set of primes dividing both $s$ and $\ell$.  Hence the period of $\mathfrak R_\ell$ divides $\operatorname{rad}(\ell)$, and its degree is at most $\ell-1$.
\end{proof}

\subsection{Odd prime indices}

\begin{theorem}[Odd prime indices]\label{thm:odd-prime-indices}
Let $p$ be an odd prime and $q-1=ps$.  Then
\[
 \mathfrak Q_p(q)=
 \begin{cases}
 \dfrac{(p-2)!\,s^{p-1}-1}{p},&p\nmid s,\\[2mm]
 \dfrac{(p-1)!\,s^{p-1}}{p^2}+p-2,&p\mid s,
 \end{cases}
 \qquad
 \mathfrak R_p(q)=
 \begin{cases}(p-2)!\,s^{p-1}-1,&p\nmid s,\\[1mm]
 \dfrac{(p-1)!\,s^{p-1}}p-1,&p\mid s,
 \end{cases}
\]
and therefore
\[
 \mathfrak A_p(q)=
 \begin{cases}
 \dfrac{p+1}{p}\bigl((p-2)!\,s^{p-1}-1\bigr),&p\nmid s,\\[2mm]
 \dfrac{p+1}{p^2}(p-1)!\,s^{p-1}+p-3,&p\mid s.
 \end{cases}
\]
\end{theorem}

\begin{proof}
Suppose first that $p\nmid s$.  Then $\mathcal A_{p,s}$ is the full affine group of $\mathbb Z/p\mathbb Z$ and $\ker\pi$ is trivial.  By \eqref{eq:exact-cardinality} the identity fixes $(p-1)!\,s^{p-1}-(p-1)$ points.  An affine map with multiplier $\delta\neq1$ has a fixed point and is not the identity, so it is not conjugate to a translation; by Theorem~\ref{thm:double-coset} it lies in no $K_\sigma$.  For a nonzero translation one has $\kappa=0$, and the fixed-point equation is $U_{j+t}=U_jU_t$.  Since $t\neq0$ generates $\mathbb Z/p\mathbb Z$, the map $j\mapsto U_j$ is a character of $\mathbb Z/p\mathbb Z$, and its inverse Fourier transform has exactly one nonzero coefficient.  Such a point corresponds to a monomial, which has index $1$ and not $p$.  Hence only the identity contributes, and Burnside's lemma gives the first formula for $\mathfrak Q_p$.

Suppose now that $p\mid s$.  Then $\delta=1$ for every element of $\Gamma^0_{p,s}$, and $\lvert\Gamma^0_{p,s}\rvert=p^2$.  If $t=0$ and $\kappa\neq0$, the fixed-point equation at $j=1$ gives $U_1=\zeta^\kappa U_1$, which is impossible.  If $\kappa=0$ and $t\neq0$, the preceding character argument again gives only monomials.  Thus only the pairs with $\kappa\neq0$ and $t\neq0$ remain.

Fix such a pair; all indices below are in $\mathbb Z/p\mathbb Z$.  Invariance of the induced permutation gives $\sigma(j+t)=\sigma(j)+\sigma(t)$, so, because $t$ generates $\mathbb Z/p\mathbb Z$, one has $\sigma(j)=wj$ for some $w\in(\mathbb Z/p\mathbb Z)^*$.  Invariance of $U$ gives
\[
 U_{j+t}=\zeta^{-\kappa j}U_jU_t,
\]
and hence, for $n\in\mathbb Z/p\mathbb Z$,
\[
 U_{nt}=U_t^n\zeta^{-\kappa tn(n-1)/2},\qquad U_t^p=1.
\]
Since $p\mid s$, the last relation gives $U_t^s=1$.  But $U_t^s=\zeta^{\sigma(t)-t}$, so $\sigma(t)=t$; hence $w=1$ and $\sigma=\operatorname{id}$.  Conversely, every choice of $U_t\in\mu_p$ determines a fixed point through the displayed recurrence, so there are exactly $p$ fixed points.

It remains to check the index.  Write $U_t=\zeta^b$ and parametrize $j=nt$.  The inverse Fourier coefficient in slot $k$ is
\[
 c_k=\frac1p\sum_{n\in\mathbb Z/p\mathbb Z}\zeta^{-\kappa tn^2/2+(b+\kappa t/2-kt)n}.
\]
The quadratic coefficient $-\kappa t/2$ is nonzero, and completing the square shows that $c_k$ is a nonzero multiple of a quadratic Gauss sum $G(A)=\sum_n\zeta^{An^2}$ with $A\neq0$.  Since $p$ is odd, the change of variables $u=x-y$, $v=x+y$ is invertible on $(\mathbb Z/p\mathbb Z)^2$ and gives
\[
 G(A)G(-A)=\sum_{u,v\in\mathbb Z/p\mathbb Z}\zeta^{Auv}=p\cdot1_{\mathbb F_q}.
\]
The characteristic of $\mathbb F_q$ is not $p$, so this product is nonzero and $G(A)\neq0$.  Therefore all $p$ Fourier coefficients are nonzero, the support gcd is $1$, and all these fixed points have exact index $p$.

The $(p-1)^2$ pairs with $\kappa,t\neq0$ thus contribute $p(p-1)^2$ fixed points in total.  Together with the identity contribution $(p-1)!\,s^{p-1}-p$ from \eqref{eq:exact-cardinality}, this gives the second formula for $\mathfrak Q_p$.  Finally, Proposition~\ref{prop:Rell} and \eqref{eq:LGW-general} give
\[
 \mathfrak R_p(q)=(p-1)!\,s^{p-1}\frac{\varphi(s)}{\varphi(ps)}-1,
\]
which is the displayed formula according as $p\nmid s$ or $p\mid s$; adding $\mathfrak Q_p$ gives $\mathfrak A_p$.
\end{proof}

\subsection{Index four}

\begin{corollary}[Index four]\label{cor:index-four-classes}
If $q-1=4s$, then
\[
 \mathfrak Q_4(q)=
 \begin{cases}
 (3s^3+s)/4,&s\text{ odd},\\
 (3s^3+4s)/8,&s\equiv2\pmod4,\\
 3s^3/8+1,&4\mid s,
 \end{cases}
 \qquad
 \mathfrak A_4(q)=
 \begin{cases}
 (15s^3-3s)/4,&s\text{ odd},\\
 (15s^3-4s)/8,&s\equiv2\pmod4,\\
 (15s^3-8s+8)/8,&4\mid s.
 \end{cases}
\]
\end{corollary}

\begin{proof}
For $\ell=4$, Theorem~\ref{thm:class-number} reduces to the following Burnside totals:
\[
\begin{array}{c|c|c|c}
 &\lvert\Gamma^0_{4,s}\rvert&\text{fixed points of the identity}&\text{fixed points of nonidentity elements}\\ \hline
 s\text{ odd}&8&6s^3-2s&4s\\
 s\equiv2\pmod4&16&6s^3-4s&12s\\
 4\mid s&16&6s^3-4s&4s+16
\end{array}.
\]
In the notation $(\kappa,t)$ of Section~\ref{sec:orbits}, the nonidentity elements with fixed points are as follows.  For $s$ odd, they are $(2,1)$ and $(2,3)$, each with $2s$ fixed points.  For $s\equiv2\pmod4$, they are $(2,2)$ and the two elements with $\kappa\equiv s/2\pmod4$ and $t\in\{1,3\}$, each with $4s$ fixed points.  For $4\mid s$, they are $(2,2)$, with $4s$ fixed points, and the four elements with $\kappa$ and $t$ odd, each with $4$ fixed points.  Division by the group order gives $\mathfrak Q_4(q)$, and Proposition~\ref{prop:Rell} then gives $\mathfrak A_4(q)$.
\end{proof}

\begin{proof}[Proof of Theorems~\textup{\ref{thm:B}} and~\textup{\ref{thm:C}}]
Theorem~\ref{thm:B}~\ref{B:formula} is Theorem~\ref{thm:class-number} together with Lemma~\ref{lem:cycle-count}, part \ref{B:nonzero} is Proposition~\ref{prop:Rell} with \eqref{eq:LGW-general}, and part \ref{B:quasi} is Theorem~\ref{thm:quasipolynomial}; note that $\gcd(\ell,s)=\gcd(\ell,\varrho)$ for $s\equiv\varrho\pmod{P_\ell}$.  Theorem~\ref{thm:C} is Theorem~\ref{thm:odd-prime-indices} together with Corollary~\ref{cor:index-four-classes}.
\end{proof}

\section{Index four: decomposition by vanishing order}\label{sec:index4}

\subsection{Notation}\label{subsec:quartic}
Throughout Sections~\ref{sec:index4}--\ref{sec:trinomial}, $q\equiv1\pmod4$, $s=(q-1)/4$, and $\zeta=\gamma^s$ has order $4$, so $\zeta^2=-1$.  We use two quartic characters.  The field-valued character is
\[
 \chi(x):=x^s\in\mu_4=\{1,\zeta,-1,-\zeta\}\qquad(x\in\mathbb F_q^*),
\]
extended by $\chi(0):=0$, and the complex-valued character $\psi$ is the one of Section~\ref{subsec:main}, normalized by $\psi(\gamma)=i$ and $\psi(0)=0$; thus $\psi(\zeta)=i^s$.  We write $\eta:=\chi^2$ for the quadratic character when working in $\mathbb F_q$, and also $\eta=\psi^2$ for its complex-valued counterpart.  For complex-valued characters $\psi_1,\psi_2$ of $\mathbb F_q^*$, extended by $0$ at $0$, we write $J(\psi_1,\psi_2):=\sum_{x\in\mathbb F_q}\psi_1(x)\psi_2(1-x)$ for their Jacobi sum; thus $J=J(\psi,\psi)$.  For a smooth projective curve $C$ over $\mathbb F_q$ we call $a_q(C)=q+1-\lvert C(\mathbb F_q)\rvert$ its Frobenius trace.

Recall the numbers $B(t)$ and $\mathfrak T(t)$ of Section~\ref{subsec:main}.  Since $u\mapsto u^t$ on $\mu_4$ depends only on $t\bmod4$, so does $B(t)$, and we write $B_c:=B(t)$ for $t\equiv c\pmod4$.  For a PP $F$ vanishing at $0$, let $v(F):=\min\{e\geqslant1:[X^e]F\neq0\}$ be its vanishing order, and let $\mathcal C_4(q,r)$ be the number of PPs of exact index $4$ vanishing at $0$ with $v(F)=r$, counted up to multiplication by nonzero constants.  Recall that $L_s=s\bigl(3+(-1)^s\bigr)$.

\subsection{Counts at a fixed residue exponent}

\begin{proposition}[Counts by vanishing order]\label{prop:support-filtration}
Let $1\leqslant a\leqslant s$ with $\gcd(a,s)=1$.  Then
\begin{align}
 \mathcal C_4(q,a)&=6s^3-\mathfrak T(a+s)-B(a+s)-B(a+2s)-L_s,\label{eq:C4-layer0}\\
 \mathcal C_4(q,a+s)&=\mathfrak T(a+s)+B(a+s),\label{eq:C4-layer1}\\
 \mathcal C_4(q,a+2s)&=B(a+2s),\label{eq:C4-layer2}
\end{align}
and no PP of exact index $4$ vanishing at $0$ with residue exponent $a$ has vanishing order $a+3s$.
\end{proposition}

\begin{proof}
Write $F=X^ah(X^s)$ with $h(T)=\sum_{k=0}^3c_kT^k$ and $A_j=h(\zeta^j)$.  After normalizing $A_0=1$, Lemma~\ref{lem:cyclo} shows that the six permutations of $\mathbb Z/4\mathbb Z$ fixing $0$, together with $s$ choices for each of $A_1,A_2,A_3$, give $6s^3$ PPs with residue exponent $a$ up to scalars.  By Proposition~\ref{prop:support-gcd}, a lower index means that the exponents $k$ with $c_k\neq0$ all have the same parity.  The two sets so obtained, with only even and with only odd such exponents, are disjoint.  Writing $h(T)=\tilde h(T^2)$ in the first case and $h(T)=T\tilde h(T^2)$ in the second reduces each count to the cyclotomic criterion for index $2$, and the two counts are $2s\mathbf 1_{\gcd(a,2s)=1}$ and $2s\mathbf 1_{\gcd(a+s,2s)=1}$.  Since $\gcd(a,s)=1$, both indicators equal $1$ if $s$ is even, whereas exactly one of them equals $1$ if $s$ is odd.  Thus
\[
 2s\bigl(\mathbf 1_{\gcd(a,2s)=1}+\mathbf 1_{\gcd(a+s,2s)=1}\bigr)=s\bigl(3+(-1)^s\bigr)=L_s.
\]

Among the points of exact index $4$, those whose least exponent $k$ with $c_k\neq0$ is $2$ have support $\{2,3\}$ and contribute $B(a+2s)$.  Those whose least such exponent is $1$ have support $\{1,2\}$ or $\{1,2,3\}$ and contribute $B(a+s)$ and $\mathfrak T(a+s)$; the support $\{1,3\}$ has index $2$, and the support $\{3\}$ gives a monomial.  This proves \eqref{eq:C4-layer1}, \eqref{eq:C4-layer2} and the last assertion, and subtracting from $6s^3-L_s$ gives \eqref{eq:C4-layer0}.  Finally, $\gcd(a+js,s)=1$ and $a+s\leqslant2s$, so the $\mathfrak T$- and $B$-terms are defined.
\end{proof}

\begin{theorem}[Sum over the vanishing orders]\label{thm:residue-collapse}
For every $a$ with $1\leqslant a\leqslant s$ and $\gcd(a,s)=1$,
\[
 \mathcal C_4(q,a)+\mathcal C_4(q,a+s)+\mathcal C_4(q,a+2s)=6s^3-L_s.
\]
Consequently, every quantity that enters these counts only through $\mathfrak T(a+s)$, $B(a+s)$ and $B(a+2s)$, in particular every quartic Jacobi sum and Frobenius trace in the formulas of Sections~\textup{\ref{sec:binomial}} and~\textup{\ref{sec:trinomial}}, cancels from this sum.
\end{theorem}

\begin{proof}
Add \eqref{eq:C4-layer0}--\eqref{eq:C4-layer2}.  The coefficients of $\mathfrak T(a+s)$, $B(a+s)$ and $B(a+2s)$ in the sum are zero.
\end{proof}

\begin{corollary}[The global count for index four]\label{cor:LGW-recovery}
The number of PPs of $\mathbb F_q$ of exact index $4$ vanishing at $0$ is
\[
 N_0(4,q)=24s^4\varphi(s)-8s^2\varphi(2s),
\]
and the number of all PPs of exact index $4$ is $N(4,q)=q\,N_0(4,q)$.  This agrees with \eqref{eq:LGW-general} for $\ell=4$, whose terms for $d=4,2,1$ are $24s^4\varphi(s)$, $-8s^2\varphi(2s)$ and $0$, and hence with Li, Guan and Wu \cite[Theorem~1.4 and its proof, pp.~841--842]{LiGuanWu2025}.
\end{corollary}

\begin{proof}
Every PP of exact index $4$ vanishing at $0$ has a unique vanishing order $a+js$ with $1\leqslant a\leqslant s$ and $j\in\{0,1,2\}$: the residue class determines $a$, Lemma~\ref{lem:cyclo} gives $\gcd(a,s)=1$, and Proposition~\ref{prop:support-filtration} excludes $j=3$.  Hence the three sets counted in Theorem~\ref{thm:residue-collapse} are disjoint and exhaustive.  Multiplying their sum by the $4s$ output scalars and summing over the $\varphi(s)$ admissible residue exponents gives $N_0(4,q)=4s\varphi(s)(6s^3-L_s)$, and $4s\varphi(s)L_s=8s^2\varphi(2s)$.  The $q$ output translations give the total count.
\end{proof}

\begin{proof}[Proof of Theorem~\textup{\ref{thm:D}~\ref{D:strata}}]
This is Proposition~\ref{prop:support-filtration} together with Theorem~\ref{thm:residue-collapse}.
\end{proof}

\section{Counting binomials}\label{sec:binomial}
In this section we evaluate the binomial count
\[
 \mathcal B_{q,r}:=B(r)+B(r+s)\qquad(\gcd(r,s)=1)
\]
that enters the trinomial count of Section~\ref{sec:trinomial}.  The letter $a$ denotes the coefficient of a binomial $X^t(X^s+a)$, with $a\in\mathbb F_q^*$, and throughout we put
\[
 \varepsilon:=\psi(\zeta)=i^s,\qquad\rho:=\psi(2).
\]

\subsection{Fourier expansion}
If $a\in\mu_4$, then $u+a=0$ for $u=-a\in\mu_4$, so $u\mapsto u^t(u+a)^s$ does not permute $\mu_4$.  For $a\in\mathbb F_q^*\setminus\mu_4$ write
\[
 \psi(a+\zeta^j)=i^{v_j(a)},\qquad v_j(a)\in\mathbb Z/4\mathbb Z.
\]
Then $\chi(a+\zeta^j)=\zeta^{v_j(a)}$, so the value of $u^t(u+a)^s$ at $u=\zeta^j$ is $\zeta^{tj+v_j(a)}$, and the induced map is a permutation of $\mu_4$ precisely when
\[
 \bigl(v_0(a),\,v_1(a)+t,\,v_2(a)+2t,\,v_3(a)+3t\bigr)
\]
is a permutation of $\mathbb Z/4\mathbb Z$.  For $m=(m_0,\dots,m_3)\in(\mathbb Z/4\mathbb Z)^4$ define
\begin{equation}\label{eq:Cm}
 C(m):=\sum_{\sigma\in S_4}i^{-\sum_{j=0}^3m_j\sigma(j)},
 \qquad
 w(m):=\sum_{j=0}^3jm_j\pmod4,
\end{equation}
where $S_4$ denotes the bijections $\{0,1,2,3\}\to\mathbb Z/4\mathbb Z$, and
\begin{equation}\label{eq:Sm}
 \mathcal S_m(q):=\sum_{a\in\mathbb F_q^*\setminus\mu_4}\prod_{j=0}^3\psi^{m_j}(a+\zeta^j).
\end{equation}
Fourier inversion on $(\mathbb Z/4\mathbb Z)^4$ gives
\begin{equation}\label{eq:B-fourier}
 B_c=\frac1{256}\sum_mC(m)\,i^{cw(m)}\mathcal S_m\qquad(c\in\mathbb Z/4\mathbb Z).
\end{equation}

\begin{lemma}[Sparse Fourier support]\label{lem:fourier-support}
For $m\in(\mathbb Z/4\mathbb Z)^4$, the coefficient $C(m)$ of \eqref{eq:Cm} satisfies
\[
 C(m)\neq0\quad\Longleftrightarrow\quad\sum_{j=0}^3m_j\equiv0\pmod4.
\]
Hence exactly $64$ tuples $m$ have $C(m)\neq0$.  Their coefficients lie in $\{\pm8,\pm24\}$, and according to the number $n(m):=\lvert\{j:m_j\neq0\}\rvert$ of nonzero coordinates they split as
\[
 64=1+18+24+21\qquad(n(m)=0,2,3,4).
\]
\end{lemma}

\begin{proof}
Adding a constant $k$ to all four values of a permutation multiplies $C(m)$ by $i^{-k\sum m_j}$, which forces the displayed congruence for a nonzero coefficient.  Checking the $4^3=64$ tuples satisfying the congruence gives the converse and the set of coefficients.  The counts by support size are elementary: $n(m)=2$ gives $\binom42\cdot3=18$ tuples, and the remaining counts follow by subtraction or by direct enumeration.
\end{proof}

Because $\mathcal S_m$ is a sum over the common punctured set $\mathbb F_q^*\setminus\mu_4$, the omitted points at which the character is unramified have to be restored.  Fix a prime number $\varpi\nmid q$.  For a tuple $m$ with $C(m)\neq0$ let
\[
 U_m:=\mathbb P^1\setminus\{-\zeta^j:m_j\neq0\},
\]
and let $\mathcal L_m$ be the rank-one $\varpi$-adic Kummer sheaf on $U_m$ whose Frobenius trace at $a$ is the character factor in \eqref{eq:Sm}.  We use Katz's construction of rank-one character sheaves from the Lang torsor and the Euler--Poincar\'e and Lefschetz trace formulas \cite[\S\S2.3 and 4.3]{Katz1988}.  We fix an embedding $\mathbb Q(i)\hookrightarrow\overline{\mathbb Q}_\varpi$ and regard $\psi$ as $\varpi$-adic.  In this section $H^i$ and $H_c^i$ denote the \'etale cohomology of the base change to $\overline{\mathbb F}_q$ with coefficients in $\overline{\mathbb Q}_\varpi$, and $\operatorname{Frob}_q$ denotes the geometric Frobenius.

\begin{proposition}[Cohomology of the Kummer sheaves]\label{prop:rank-collapse}
For $m\neq0$ with $C(m)\neq0$,
\[
 \dim H_c^1(U_m,\mathcal L_m)=n(m)-2.
\]
Thus tuples with $n(m)=2,3,4$ give cohomological dimensions $0,1,2$; the points omitted from the summation domain of $\mathcal S_m$ are treated below.
\end{proposition}

\begin{proof}
The order of the character divides $4$, so $\mathcal L_m$ is tame.  Its local monodromy is nontrivial at exactly the $n(m)$ finite punctures, while $\sum m_j\equiv0\pmod4$ makes $\mathcal L_m$ unramified at infinity.  Hence the tame Euler--Poincar\'e formula gives $\chi_c(U_m,\mathcal L_m)=2-n(m)$.  Since $U_m$ is affine, $H_c^0=0$, and since $\mathcal L_m$ is a nontrivial rank-one local system, it has no geometric coinvariants, so $H_c^2=0$ \cite[2.0.5--2.0.6]{Katz1988}.  This proves the formula.
\end{proof}

\subsection{Support four and Frobenius traces}\label{subsec:support4}
Let $m$ be a tuple with $C(m)\neq0$ and $n(m)=4$.  Then $U_m(\mathbb F_q)$ contains, besides the points of $\mathbb F_q^*\setminus\mu_4$, the points $a=0$ and $a=\infty$, which are not in the summation domain of \eqref{eq:Sm}.  The sheaf $\mathcal L_m$ is unramified at both, with Frobenius traces $\prod_j\psi(\zeta^j)^{m_j}=\varepsilon^{w(m)}$ and $1$, respectively.  By the Lefschetz trace formula and Proposition~\ref{prop:rank-collapse},
\begin{equation}\label{eq:completed-trace}
 \mathcal S_m+1+\varepsilon^{w(m)}=-\operatorname{Tr}\bigl(\operatorname{Frob}_q\mid H_c^1(U_m,\mathcal L_m)\bigr).
\end{equation}

Suppose further that some $m_j$ is odd, that is, $\gcd(4,m_0,m_1,m_2,m_3)=1$; then the local monodromy of $\mathcal L_m$ at $-\zeta^j$ has order $4$, so the geometric monodromy of $\mathcal L_m$ has order $4$.  Let $Y_m$ be the smooth projective model of the cyclic cover $z^4=\prod_j(a+\zeta^j)^{m_j}$, and let $\tau\colon(a,z)\mapsto(a,\zeta z)$ be its automorphism of order $4$.  The finite cyclic-cover decomposition \cite[Lemma~7.4 and the discussion following it]{Katz2009Dwork} identifies $H_c^1(U_m,\mathcal L_m)$ with one faithful $\tau$-eigenspace $V_i\subset H^1(Y_m)$, and the conjugate Kummer character gives $V_{-i}$.  Indeed, every finite puncture has nontrivial local monodromy, so no boundary invariant stalk occurs in either faithful character; this is the decomposition of the cohomology of the cyclic cover into eigenspaces of $\tau$.  By Proposition~\ref{prop:rank-collapse}, $\dim V_i=\dim V_{-i}=2$.

\begin{lemma}[Frobenius traces through a quotient]\label{lem:trace-transfer}
Let $m$ be a tuple with $C(m)\neq0$ and $n(m)=4$ whose character has order $4$, and let $Y_m$ and $\tau$ be as above.  Suppose that an involution $\iota$ of $Y_m$ defined over $\mathbb F_q$ satisfies $\iota\tau\iota^{-1}=\tau^{-1}$, and that the smooth projective quotient $E:=Y_m/\langle\iota\rangle$ satisfies $\dim H^1(E)=2$.  Then
\[
 \mathcal S_m+1+\varepsilon^{w(m)}=-a_q(E).
\]
\end{lemma}

\begin{proof}
The involution $\iota$ exchanges $V_i$ and $V_{-i}$, so $v\mapsto v+\iota(v)$ embeds $V_i$ as a $2$-dimensional subspace of $H^1(Y_m)^{\langle\iota\rangle}$.  The quotient theorem gives a Frobenius-compatible isomorphism $H^1(E)\cong H^1(Y_m)^{\langle\iota\rangle}$ \cite[Theorem~3.1(2)--(4)]{DokchitserDokchitser2018}.  Both sides have dimension $2$, so the image of $v\mapsto v+\iota(v)$ is the whole invariant subspace.  Since $\iota$ and $\tau$ are defined over $\mathbb F_q$, Frobenius commutes with them, and its trace on this subspace equals its trace on $V_i$.  As $a_q(E)=\operatorname{Tr}(\operatorname{Frob}_q\mid H^1(E))$, the assertion follows from \eqref{eq:completed-trace}.
\end{proof}

\begin{lemma}[The tuple $1133$]\label{lem:S1133}
For every $q\equiv1\pmod4$,
\begin{equation}\label{eq:S1133}
 \mathcal S_{1133}+2=-a_q(E_{128}).
\end{equation}
\end{lemma}

\begin{proof}
For $m=(1,1,3,3)$,
\[
 \mathcal S_{1133}=\sum_{a\in\mathbb F_q^*\setminus\mu_4}\psi\!\left(\frac{(a+1)(a+\zeta)}{(a-1)(a-\zeta)}\right).
\]
Here $z^4=(a+1)(a+\zeta)(a-1)^3(a-\zeta)^3$.  Put $t:=z/\bigl((a-1)(a-\zeta)\bigr)$, so that
\[
 t^4=\frac{(a+1)(a+\zeta)}{(a-1)(a-\zeta)}
\]
and $\tau(t)=\zeta t$.  Thus $a$ is a root of
\[
 (t^4-1)a^2-(1+\zeta)(t^4+1)a+\zeta(t^4-1)=0,
\]
whose discriminant is $-2\zeta(t^8-6t^4+1)$.  Hence the function field of $Y_m$ is generated by $t$ and $z':=2(t^4-1)a-(1+\zeta)(t^4+1)$, which satisfy
\[
 z'^2=-2\zeta(t^8-6t^4+1),
\]
and $\tau(t,z')=(\zeta t,z')$.  The involution $\iota\colon a\mapsto-a$ sends $t^4$ to $t^{-4}$; in these coordinates it is $\iota\colon(t,z')\mapsto(t^{-1},z'/t^4)$.  Put $x=t+t^{-1}$ and $y=z'/t^2$; the quotient is
\[
 y^2=-2\zeta(x^4-4x^2-4).
\]
Since $(1-\zeta)^2=-2\zeta$, take $c=1-\zeta$ and put $v=y/c$.  From $v^2=x^4-4x^2-4$, the substitutions
\[
 U=\frac{v+2\zeta}{x^2},\qquad W=x(U^2-1),\qquad X=4\zeta U-4,\qquad Y=4\zeta W
\]
give the curve $E_{128}$.  The involution $\iota$ satisfies $\iota\tau\iota^{-1}=\tau^{-1}$ for the automorphism $\tau$ above, and the smooth projective quotient has two-dimensional $H^1$, so Lemma~\ref{lem:trace-transfer} applies.  Here $w(1133)=0$, so both local terms equal $1$.
\end{proof}

Here $j(E)$ denotes the $j$-invariant of the smooth projective curve defined by a nonsingular cubic of the above form, with its point at infinity as origin, in the convention of \cite[Chap.~III, \S1]{Silverman2009}.  Direct substitution gives $j(E_{128})=128$, $j(E_{1728})=1728$ and $j(E_{8000})=8000$.  For these curves over $\mathbb Q$, $\operatorname{End}_{\mathbb C}(E)$ denotes the endomorphism ring of the group law over $\mathbb C$, and $E$ has \emph{complex multiplication} (CM) when $\operatorname{End}_{\mathbb C}(E)\neq\mathbb Z$.  Table~1 of Bilu, Luca and Pizarro-Madariaga \cite{BiluLucaPizarro2016} lists the rational $j$-invariants arising from CM.  Since $128$ does not occur in it, $E_{128}$ has no CM, while the same table shows that $\operatorname{End}_{\mathbb C}(E_{8000})$ has discriminant $-8$.  The curve $E_{1728}$ has CM, because $(X,Y)\mapsto(-X,iY)$ is an automorphism of order $4$.

\begin{lemma}[Further completed support-four traces]\label{lem:support4-traces}
For every $q\equiv1\pmod4$,
\begin{equation}\label{eq:2222}
 \mathcal S_{2222}+2=-a_q(E_{1728}).
\end{equation}
If $s$ is even, choose $\zeta_8\in\mathbb F_q$ with $\zeta_8^2=\zeta$.  Then
\begin{align}
 \operatorname{Re}\mathcal S_{1111}+2&=-a_q(E_{1728}),\label{eq:1111-even}\\
 \mathcal S_{1313}+2&=-a_q(E_{\zeta_8}),\notag\\
 \mathcal S_{1232}+2&=-a_q(E_{8000}),\notag
\end{align}
where
\[
 E_{\zeta_8}\colon Y^2=X(X^2-8\zeta_8 X+8\zeta_8^2),
\]
and $a_q(E_{\zeta_8})=i^s\,a_q(E_{8000})$.
\end{lemma}

\begin{proof}
For \eqref{eq:2222}, completion gives the quadratic cover $y^2=x^4-1$.  The rational change of variables
\[
 X=2\,\frac{x+1}{x-1},\qquad Y=\frac{4y}{(x-1)^2}
\]
gives $Y^2=X^3+4X$, and the two completed points account for the term $2$.

Assume now that $s$ is even, so that $q\equiv1\pmod8$, $\psi(-1)=\psi(4)=1$ and $\eta(2)=1$, and choose $\beta_2\in\mathbb F_q$ with $\beta_2^2=2$.  We treat $1111$ directly, because its quotient curve records the real part rather than the whole faithful Kummer trace.  Since $\prod_{j=0}^3(a+\zeta^j)=a^4-1$, completion gives
\[
 \mathcal S_{1111}+2=1+\sum_{a\in\mathbb F_q}\psi(1-a^4).
\]
The number of fourth roots of $x\in\mathbb F_q$ is $1+\psi(x)+\eta(x)+\overline{\psi(x)}$; hence
\[
 \mathcal S_{1111}+2=J(\psi,\psi)+J(\eta,\psi).
\]
The change of variables $y=2x-1$ in $J(\psi,\psi)$ gives $J(\psi,\psi)=\overline{\psi(4)}J(\eta,\psi)=J(\eta,\psi)$, so $\operatorname{Re}\mathcal S_{1111}+2=2\operatorname{Re}J$.  The rational change of variables
\[
 U=\frac{y+1}{u^2},\qquad W=u(U^2-1),\qquad X=2\zeta U,\qquad Y=\beta_2(1-\zeta)W
\]
sends the quartic $y^2=u^4+1$ to $Y^2=X^3+4X$.  The quartic has two points at infinity, and the same fourth-root count gives
\[
 -a_q(E_{1728})=1+\sum_{u\in\mathbb F_q}\eta(1+u^4)=J(\psi,\eta)+J(\overline\psi,\eta)=2\operatorname{Re}J.
\]
This proves \eqref{eq:1111-even}.

For $1313$, let $u$ denote the fourth-root coordinate of the cover, and put $t=u/(a^2+1)$ and $z=2a/(a^2+1)$.  Then $z^2=1-t^8$, and the automorphism of order $4$ is $\tau(t,z)=(\zeta t,z)$.  The involution
\[
 \iota(t,z)=\bigl(\zeta_8/t,\ \zeta z/t^4\bigr)
\]
satisfies $\iota\tau\iota^{-1}=\tau^{-1}$.  Its invariant coordinates
\[
 x=t+\frac{\zeta_8}{t},\qquad y=\frac{z}{t^2}
\]
give the quotient
\[
 y^2=-x^4+4\zeta_8 x^2-2\zeta.
\]
Put $c=1-\zeta$.  The rational point $(0,c)$ gives the explicit conversion
\[
 U=\frac{y+c}{x^2},\qquad W=x(U^2+1),\qquad X=2(cU+2\zeta_8),\qquad Y=2cW,
\]
for which $Y^2=X(X^2-8\zeta_8 X+8\zeta_8^2)$.  Thus the quotient is $E_{\zeta_8}$.  Since $w(1313)=2$ and $\varepsilon^2=1$, Lemma~\ref{lem:trace-transfer} gives $\mathcal S_{1313}+2=-a_q(E_{\zeta_8})$.

For $1232$, after dividing the fourth-root coordinate by $a-1$, the cyclic-cover equation becomes
\[
 t^4=\frac{(a+1)(a^2+1)^2}{a-1}.
\]
Putting $\xi=t^2/(a^2+1)$ and $V=t(\xi^2-1)$ gives $V^2=2\xi(\xi^4+1)$, and the automorphism becomes $\tau(\xi,V)=(-\xi,\zeta V)$.  The involution
\[
 \iota(\xi,V)=\bigl(\xi^{-1},\ V/\xi^3\bigr)
\]
is defined over $\mathbb F_q$ and satisfies $\iota\tau\iota^{-1}=\tau^{-1}$.  With the invariant coordinates
\[
 x=\xi+\xi^{-1},\qquad y=\frac{V(\xi+1)}{\beta_2\xi^2},
\]
the quotient is $y^2=(x+2)(x^2-2)=x^3+2x^2-2x-4$.  Again $w(1232)=2$ and $\varepsilon^2=1$, so Lemma~\ref{lem:trace-transfer} gives $\mathcal S_{1232}+2=-a_q(E_{8000})$.

Finally, $E_{\zeta_8}$ has $c_4=640\zeta_8^2$ and $c_6=14336\zeta_8^3$, while $E_{8000}$ has $c_4=160$ and $c_6=1792$.  Thus $E_{\zeta_8}$ is the quadratic twist of $E_{8000}$ by $2\zeta_8$.  Since $q\equiv1\pmod8$, $\eta(2)=1$ and $\eta(\zeta_8)=\zeta_8^{2s}=\zeta^s=i^s$, so the quadratic character of the twist parameter is $i^s$.
\end{proof}

\subsection{The terms with \texorpdfstring{$w(m)=0$}{w(m)=0}}
For $\nu\in\mathbb Z/4\mathbb Z$ put
\[
 \mathcal W_\nu:=\sum_{w(m)=\nu}C(m)\mathcal S_m.
\]
Then \eqref{eq:B-fourier} becomes
\[
 B_c=\frac1{256}\sum_{\nu=0}^3i^{c\nu}\mathcal W_\nu.
\]
Two symmetries reduce the number of sums to be evaluated.  The substitution $a\mapsto\zeta a$ preserves $\mathbb F_q^*\setminus\mu_4$ and, since $\sum_jm_j\equiv0\pmod4$, shows that $\mathcal S_m$, $C(m)$ and $w(m)$ are unchanged when the coordinates of $m$ are permuted cyclically; we call the cyclic permutations of the coordinates \emph{rotations}.  Complex conjugation gives $\mathcal S_{-m}=\overline{\mathcal S_m}$.  The $16$ tuples with $C(m)\neq0$ and $w(m)=0$ are the zero tuple, two tuples with $n(m)=2$, two mutually conjugate rotation orbits of tuples with $n(m)=3$, the four rotations of $1133$, and the tuple $2222$.  Rotation and conjugation therefore give
\[
 \mathcal W_0=24(q-5)-16\mathcal S_{0202}+64\operatorname{Re}\mathcal S_{0121}+32\mathcal S_{1133}+24\mathcal S_{2222}.
\]
Since $\mathcal S_{0202}=-2-2\eta(2)$, equations \eqref{eq:S1133} and \eqref{eq:2222} give
\begin{equation}\label{eq:W0-general}
 \frac{\mathcal W_0}{256}=\frac{3q-25+4\eta(2)}{32}+\frac14\operatorname{Re}\mathcal S_{0121}-\frac18a_q(E_{128})-\frac3{32}a_q(E_{1728}).
\end{equation}
This identity is specialized below to the two residue classes of $q$ modulo $8$.

\subsection{The case \texorpdfstring{$q\equiv1\pmod8$}{q congruent to 1 mod 8}}\label{subsec:even-boundary}
Assume that $q\equiv1\pmod8$.  Then $s$ is even, $\varepsilon=i^s\in\{\pm1\}$ and $\rho=\psi(2)\in\{\pm1\}$; recall that $J=J(\psi,\psi)=\sum_{x\in\mathbb F_q}\psi(x)\psi(1-x)$.

\begin{lemma}[The support-three term]\label{lem:S0121-even}
If $q\equiv1\pmod8$, then
\[
 \mathcal S_{0121}=\rho(J+1)-2(1+\rho).
\]
\end{lemma}

\begin{proof}
Let
\[
 \Xi:=\sum_{a\in\mathbb F_q}\psi(a^2+1)\eta(a-1).
\]
The punctured sum satisfies $\mathcal S_{0121}=\Xi-1-\rho$.  For $a\neq-1$ substitute
\[
 t=\frac{a-1}{a+1}.
\]
Since
\[
 a-1=\frac{2t}{1-t},\qquad a^2+1=\frac{2(1+t^2)}{(1-t)^2}
\]
and $\eta(2)=1$, completion at the two endpoints gives
\[
 \Xi=\rho\mathcal K-1+\rho,\qquad \mathcal K:=\sum_t\psi(1+t^2)\eta(t).
\]
Using $\eta(t)=\psi(t^2)$ and the substitution $x=-t^2$,
\[
 \mathcal K=J(\psi,\psi)+J(\psi^{-1},\psi)=J-1,
\]
because $\psi(-1)=1$.  Substituting gives the lemma.
\end{proof}

The tuples with $n(m)=2$ and $w(m)=2$ are represented by $0022$ and $0103$.  Direct evaluation of the corresponding quadratic and linear-fractional character sums gives
\[
 \mathcal S_{0022}=-2-2\varepsilon\rho,\qquad \mathcal S_{0103}=-2-2\varepsilon,
\]
so their total contribution to $\mathcal W_2$ is $128+64\varepsilon(\rho+1)$.  The tuples with $n(m)=4$ and $w(m)=2$ split into the pair $1111,3333$, the pair $1313,3131$, and the four rotations of $1232$; their contribution is
\[
 -48\operatorname{Re}\mathcal S_{1111}+16\mathcal S_{1313}-32\mathcal S_{1232}=48a_q(E_{1728})+(32-16\varepsilon)a_q(E_{8000})+128
\]
by Lemma~\ref{lem:support4-traces}.  Therefore
\begin{equation}\label{eq:W2-even}
 \frac{\mathcal W_2}{128}=2+\frac{\varepsilon(\rho+1)}2+\frac38a_q(E_{1728})+\frac{2-\varepsilon}{8}a_q(E_{8000}).
\end{equation}
Inserting Lemma~\ref{lem:S0121-even} into \eqref{eq:W0-general} gives
\begin{equation}\label{eq:W0-even}
 \frac{\mathcal W_0}{128}=\frac{3q-37}{16}-\rho+\frac12\operatorname{Re}\bigl(\rho(J+1)\bigr)-\frac14a_q(E_{128})-\frac3{16}a_q(E_{1728}).
\end{equation}

If $\gcd(r,s)=1$ and $s$ is even, then $r$ is odd.  If $s\equiv0\pmod4$, then $u^s=1$ for $u\in\mu_4$, and the substitution $a\mapsto a^{-1}$, which turns $u^t(u+a^{-1})^s$ into $a^{-s}v^{-t}(a+v)^s$ with $v=u^{-1}$, gives $B_1=B_3$; if $s\equiv2\pmod4$, the exponents $r$ and $r+s$ occupy the residues $1$ and $3$.  In either case
\begin{equation}\label{eq:even-boundary-projection}
 B(r)+B(r+s)=B_1+B_3=\frac{\mathcal W_0-\mathcal W_2}{128}.
\end{equation}

\begin{theorem}[Binomial counts for $q\equiv1\pmod8$]\label{thm:even-boundary}
Assume that $q\equiv1\pmod8$, $s=(q-1)/4$ and $\gcd(r,s)=1$.  Then
\begin{align*}
 \mathcal B_{q,r}=B(r)+B(r+s)
 &=\frac{3q-69}{16}+\frac{\rho}{2}(\operatorname{Re}J-1)-\frac{\varepsilon(\rho+1)}2\\
 &\quad-\frac14a_q(E_{128})-\frac9{16}a_q(E_{1728})-\frac{2-\varepsilon}{8}a_q(E_{8000}).
\end{align*}
In particular, this boundary sum does not depend on $r$, subject to $\gcd(r,s)=1$.
\end{theorem}

\begin{proof}
Subtract \eqref{eq:W2-even} from \eqref{eq:W0-even} and use \eqref{eq:even-boundary-projection}.
\end{proof}

\subsection{The case \texorpdfstring{$q\equiv5\pmod8$}{q congruent to 5 mod 8}}\label{subsec:odd-boundary}
Assume now that $q\equiv5\pmod8$.  Then $s$ is odd and $\varepsilon=i^s\in\{\pm i\}$; put $\vartheta:=\psi(1-\zeta)$, and again $J=J(\psi,\psi)$.
For any residue $c\bmod4$, Fourier projection gives
\[
 B_c+B_{c+s}=\frac1{256}\sum_{\nu=0}^3i^{c\nu}(1+i^{s\nu})\mathcal W_\nu.
\]
Since $s$ is odd, the factor at $\nu=2$ vanishes.  Conjugation gives $\mathcal W_3=\overline{\mathcal W_1}$, and the substitution $a\mapsto a^{-1}$ gives $\mathcal W_1=\varepsilon\overline{\mathcal W_1}$; hence
\begin{equation}\label{eq:odd-projector}
 B_c+B_{c+s}=\frac{\mathcal W_0}{128}+\frac1{128}\operatorname{Re}\bigl(i^c(1+\varepsilon)\mathcal W_1\bigr).
\end{equation}

The $16$ tuples with $C(m)\neq0$ and $w(m)=1$ form four rotation orbits, with representatives
\[
 \mathcal P:=\mathcal S_{0031},\qquad \mathcal Q:=\mathcal S_{0112},\qquad \mathcal R:=\mathcal S_{0233},\qquad\Theta:=\mathcal S_{1322},
\]
and Fourier coefficients $-8,+8,+8,-8$.  Thus
\begin{equation}\label{eq:W1-PQRT}
 \mathcal W_1=32(-\mathcal P+\mathcal Q+\mathcal R-\Theta).
\end{equation}
Here $\psi(-1)=\eta(2)=-1$, $\varepsilon^2=\rho^2=-1$, and $\vartheta^2=-\rho\varepsilon$.  The substitutions are
\[
\begin{array}{c|c|c}
 &x&\text{transformed summand; omitted finite values of }x\\ \hline
 \mathcal P&(a-\zeta)/(a-1)&\psi(x);\quad 0,\,1,\,\zeta,\,1+\zeta,\,(1+\zeta)/2\\
 \mathcal Q&-\zeta(a+\zeta)/(a-\zeta)&\rho\varepsilon\vartheta\,\psi(x)\psi(x-1);\quad 0,\,1,\,-1,\,\zeta,\,-\zeta\\
 \mathcal R&-\zeta(a+\zeta)/(a-\zeta)&-\rho\vartheta^{-1}\eta(x)\psi^{-1}(1-x);\quad 0,\,1,\,-1,\,\zeta,\,-\zeta
\end{array}.
\]
The corresponding complete sums are $0$, $-J$ and $-\overline J$, respectively; the last identity follows from $x\mapsto x^{-1}$.  Subtracting the omitted terms and using $\vartheta^2=-\rho\varepsilon$, in the forms $-\rho\varepsilon\vartheta=\vartheta^{-1}$ and $\rho\vartheta^{-1}=\varepsilon\vartheta$, gives
\begin{equation}\label{eq:PQR}
\begin{aligned}
 \mathcal P&=-1-\varepsilon-\vartheta^{-1}-\varepsilon\vartheta,\\
 \mathcal Q&=\vartheta^{-1}J-1-\varepsilon+\varepsilon\vartheta,\\
 \mathcal R&=-1+\varepsilon\vartheta\overline J-\varepsilon+\vartheta^{-1}.
\end{aligned}
\end{equation}

The remaining orbit is represented by $1322$.  If $z$ is the cover coordinate, then
\[
 z^4=\frac{(a+1)(a-1)^2(a-\zeta)^2}{a+\zeta}.
\]
Putting
\[
 x:=\frac{z^2}{(a-1)(a-\zeta)},\qquad y:=z(x^2-1)
\]
gives the curve
\[
 H_q\colon\quad y^2=x\,[2-(1+\zeta)x^2]\,[1+\zeta-2\zeta x^2].
\]
Because $\zeta^2=-1$ and the characteristic is odd, $1+\zeta\neq0$, so each quadratic factor has two distinct nonzero roots.  The two factors have no common root, since a common root would give $4\zeta=(1+\zeta)^2=2\zeta$, which is impossible in odd characteristic.  Neither factor vanishes at $x=0$, and the leading coefficient $2\zeta(1+\zeta)$ is nonzero.  Thus the polynomial of degree $5$ on the right is squarefree, and its smooth projective model $H_q$ has genus $2$.

The birational change of variables is explicit in both directions:
\[
 x^2=\frac{a+1}{a+\zeta},\qquad a=\frac{1-\zeta x^2}{x^2-1},\qquad z=\frac{y}{x^2-1}.
\]
In these coordinates $\tau(x,y)=(-x,\zeta y)$, so $\tau^2(x,y)=(x,-y)$ is the hyperelliptic involution.  Its quotient is $\mathbb P^1_x$, hence $H^1(H_q)^{\langle\tau^2\rangle}=0$; in particular, the eigenspaces of $\tau$ for the eigenvalues $1$ and $-1$ vanish.  The two faithful eigenspaces $V_i$ and $V_{-i}$ each have dimension $2$ by Proposition~\ref{prop:rank-collapse}, since the character and its conjugate have four nonzero coordinates.  Thus $H^1(H_q)=V_i\oplus V_{-i}$ and $\dim H^1(H_q)=4$.  Since $w(1322)=1$, \eqref{eq:completed-trace} gives
\[
 \widetilde\Theta:=\Theta+1+\varepsilon=-\operatorname{Tr}\bigl(\operatorname{Frob}_q\mid V_i\bigr).
\]
The conjugate character gives the trace on $V_{-i}$, so the sum of the two faithful eigenspace traces is the trace on all of $H^1(H_q)$:
\[
 \widetilde\Theta+\overline{\widetilde\Theta}=-a_q(H_q).
\]
Inversion $a\mapsto a^{-1}$ sends $1322$ to $1223$ with multiplier $\varepsilon$, and a rotation sends $1223$ to $3122=\overline{1322}$.  These operations preserve the common punctured summation domain, so $\Theta=\varepsilon\overline\Theta$.  Because $\varepsilon\overline{(1+\varepsilon)}=1+\varepsilon$ when $\varepsilon^2=-1$, the completed trace satisfies $\widetilde\Theta=\varepsilon\overline{\widetilde\Theta}$ as well.  Hence
\begin{equation}\label{eq:T-Hq}
 \Theta=-(1+\varepsilon)\left(\frac{a_q(H_q)}2+1\right).
\end{equation}
Substituting \eqref{eq:PQR}--\eqref{eq:T-Hq} into \eqref{eq:W1-PQRT}, we may write
\begin{equation}\label{eq:Lambda}
 \mathcal W_1=32(1+\varepsilon)\Lambda_q,
\end{equation}
where
\begin{equation}\label{eq:Lambda-explicit}
 \Lambda_q:=\frac{\vartheta^{-1}(J+2)+\varepsilon\vartheta(\overline J+2)}{1+\varepsilon}+\frac{a_q(H_q)}2.
\end{equation}
The number $\Lambda_q$ is an integer.  Indeed, each of $\mathcal P,\mathcal Q,\mathcal R,\Theta$ lies in $\mathbb Z[i]$, so $Z:=\mathcal W_1/32=-\mathcal P+\mathcal Q+\mathcal R-\Theta\in\mathbb Z[i]$.  The symmetry $Z=\varepsilon\overline Z$ with $\varepsilon=\pm i$ forces $Z\in(1+\varepsilon)\mathbb Z$, and hence $\Lambda_q=Z/(1+\varepsilon)\in\mathbb Z$ by \eqref{eq:Lambda}.

The support-three term with $w(m)=0$ is now
\[
 \mathcal S_{0121}=\psi(2)(J+1)-2.
\]
Since $\eta(2)=-1$, \eqref{eq:W0-general} becomes
\begin{equation}\label{eq:Mq}
 M_q:=\frac{\mathcal W_0}{128}=\frac{3q-45}{16}+\frac12\operatorname{Re}\bigl(\psi(2)(J+1)\bigr)-\frac14a_q(E_{128})-\frac3{16}a_q(E_{1728}).
\end{equation}

\begin{theorem}[Binomial counts for $q\equiv5\pmod8$]\label{thm:odd-boundary}
Assume that $q\equiv5\pmod8$, $s=(q-1)/4$ and $\gcd(r,s)=1$.  Then
\[
 \mathcal B_{q,r}=B(r)+B(r+s)=M_q+\frac{\Lambda_q}{2}\operatorname{Re}(i^{r+s}),
\]
with $M_q$ and $\Lambda_q$ given by \eqref{eq:Mq} and \eqref{eq:Lambda-explicit}.
\end{theorem}

\begin{proof}
Use \eqref{eq:odd-projector} with $c\equiv r\pmod4$, together with \eqref{eq:Lambda}.  Since $(1+\varepsilon)^2=2\varepsilon$ and $\varepsilon=i^s$, the real part simplifies to the displayed term.
\end{proof}

\begin{remark}[Individual binomial counts]\label{rem:individual-B}
Replacing $b$ by $b^{-1}$ and putting $v=u^{-1}$, one has $u^t(u+b^{-1})^s=b^{-s}v^{-t-s}(b+v)^s$, so $B(t)=B(-t-s)$ for every integer $t$.  If $q\equiv5\pmod8$, then $r$ in Theorem~\ref{thm:odd-boundary} may lie in any residue class modulo $4$, and this symmetry gives $B_0=B_3$ and $B_1=B_2$ for $s\equiv1\pmod4$, and $B_0=B_1$ and $B_2=B_3$ for $s\equiv3\pmod4$; hence the sums $B_c+B_{c+s}$ determine every $B_c$.  If $q\equiv1\pmod{16}$, then $4\mid s$ and \eqref{eq:even-boundary-projection} gives $B_1=B_3=(\mathcal W_0-\mathcal W_2)/256$.  If $q\equiv9\pmod{16}$, the symmetry fixes the residues $1$ and $3$, and Theorem~\ref{thm:even-boundary} determines only $B_1+B_3$.  The difference $B_1-B_3$ is given by \eqref{eq:B-fourier}, but it is not evaluated here, and it need not vanish: for $q=9$ one has $B_1=4$ and $B_3=0$, and for $q=73$ and $q=89$ one has $(B_1,B_3)=(4,8)$ and $(0,12)$.  Accordingly, for $q\equiv9\pmod{16}$ the counts in Theorem~\ref{thm:D}~\ref{D:strata} involve the individual numbers $B(a+s)$ and $B(a+2s)$, which are not reduced to Jacobi sums and Frobenius traces here.
\end{remark}

\section{Counting trinomials}\label{sec:trinomial}
In this section $1\leqslant r\leqslant2s$ and $\gcd(r,s)=1$.  We count the projective classes $[a:b:c]$ with $abc\neq0$ for which
\[
 X^r(a+bX^s+cX^{2s})
\]
permutes $\mathbb F_q$; their number is $\mathfrak T(r)$, and $r$ is the vanishing order of these trinomials.  We first count all classes with $b\neq0$, allowing $a=0$ or $c=0$; the classes with $a=0$ or $c=0$, which correspond to binomials and to one monomial, are removed in Section~\ref{subsec:trinomial-count}.  By Lemma~\ref{lem:cyclo} with $\ell=4$, for $e\geqslant1$ and $h\in\mathbb F_q[T]$ the polynomial $X^eh(X^s)$ permutes $\mathbb F_q$ if and only if $\gcd(e,s)=1$ and $u\mapsto u^eh(u)^s$ is a permutation of $\mu_4$.

\subsection{A criterion on \texorpdfstring{$\mu_4$}{mu4}}
For the trinomial family write
\begin{equation}\label{eq:Aj}
 h(T):=a+bT+cT^2,\qquad A_j=h(\zeta^j)=a+b\zeta^j+c\zeta^{2j}.
\end{equation}
The induced values are then
\begin{equation}\label{eq:yj}
 y_j:=\zeta^{rj}\chi(A_j)\in\mu_4\cup\{0\},
\end{equation}
and $y_j=0$ exactly when $A_j=0$.

\begin{lemma}[Sum and product on $\mu_4$]\label{lem:sum-product}
Let $y_0,y_1,y_2,y_3\in\mu_4$, in odd characteristic.  Then $\{y_0,y_1,y_2,y_3\}=\mu_4$ if and only if
\[
 \sum_{j=0}^3y_j=0\qquad\text{and}\qquad\prod_{j=0}^3y_j=-1.
\]
\end{lemma}

\begin{proof}
The forward implication is immediate.  Conversely, let $e_k$ be the elementary symmetric functions and $p_k=\sum_jy_j^k$.  We have $e_1=0$, $e_4=-1$ and $p_4=4$.  Newton's identities give $p_2=-2e_2$ and
\[
 4=p_4=-e_2p_2-4e_4=2e_2^2+4,
\]
so $e_2=0$.  Since $y_j^{-1}=y_j^3$,
\[
 e_3=e_4\sum_jy_j^{-1}=-p_3,
\]
whereas Newton's identities give $p_3=3e_3$.  Thus $4p_3=0$, whence $e_3=p_3=0$ in odd characteristic.  Therefore $\prod_j(T-y_j)=T^4-1$.
\end{proof}

By \eqref{eq:yj} and Lemma~\ref{lem:sum-product}, $X^r(a+bX^s+cX^{2s})$ is a PP if and only if all $A_j$ are nonzero and
\[
 \sum_{j=0}^3\zeta^{rj}\chi(A_j)=0\qquad\text{and}\qquad
 \chi\bigl(((a+c)^2-b^2)((a-c)^2+b^2)\bigr)=(-1)^{r+1};
\]
the second condition follows from
\[
 \prod_{j=0}^3A_j=\bigl((a+c)^2-b^2\bigr)\bigl((a-c)^2+b^2\bigr)\qquad\text{and}\qquad\sum_{j=0}^3j=6.
\]

\subsection{Opposite pairs}\label{subsec:opposite}
The four elements of $\mu_4$ form the two opposite pairs $\{\pm1\}$ and $\{\pm\zeta\}$.  Hence, if $(y_0,y_1,y_2,y_3)$ is a permutation of $\mu_4$, the set of positions $j$ with $y_j\in\{\pm1\}$ and its complement form one of the three partitions $\{\{0,2\},\{1,3\}\}$, $\{\{0,1\},\{2,3\}\}$ and $\{\{0,3\},\{1,2\}\}$.  We say that $y$ is of \emph{diagonal type} if the pairs are $\{0,2\},\{1,3\}$, of \emph{edge type} $+$ if they are $\{0,1\},\{2,3\}$, and of \emph{edge type} $-$ if they are $\{0,3\},\{1,2\}$; every permutation has exactly one type.  Under $h(T)\mapsto h(\zeta T)$ the induced values become $y'_j=\zeta^{-r}y_{j+1}$.  Thus a cyclic relabelling of the four positions, together with a common nonzero scalar, preserves the opposite-pair relation; it preserves the diagonal type and exchanges the two edge types.  Hence the two edge types occur equally often, and the count before the binomial boundary is removed is
\begin{equation}\label{eq:D2E}
 N_{\mathrm{pre}}(r)=D_r+2E_r,
\end{equation}
where $D_r$ and $E_r$ denote the numbers of projective classes $[a:b:c]$ with $b\neq0$, and $a,c$ arbitrary, of diagonal type and of edge type $+$, respectively.  We normalize $b=1$.  The values \eqref{eq:Aj} satisfy the Fourier relation
\begin{equation}\label{eq:fourier-relation}
 A_0+\zeta A_1-A_2-\zeta A_3=0,
\end{equation}
which expresses the vanishing of the coefficient of $T^3$.

The \emph{quartic cyclotomic classes} and \emph{quartic cyclotomic numbers} are
\[
 C_j:=\{x\in\mathbb F_q^*:\chi(x)=\zeta^j\},\qquad
 (i,j)_4:=\lvert\{x\in C_i:x+1\in C_j\}\rvert.
\]

\subsection{The diagonal type}
Put
\[
 \varepsilon_D:=(-1)^{r+1},\qquad
 k_D:=\begin{cases}0,&r\text{ odd},\\2,&r\text{ even},\end{cases}
\]
so that $\varepsilon_D=\zeta^{k_D}$, and let
\[
 U:=A_0/A_2,\qquad V:=A_1/A_3.
\]
The conditions for the diagonal type become
\[
 \chi(U)=\chi(V)=\varepsilon_D,\qquad\eta(A_2/A_3)=\varepsilon_D.
\]
By \eqref{eq:fourier-relation},
\begin{equation}\label{eq:D-ratio}
 \frac{A_2}{A_3}=-\zeta\,\frac{V-1}{U-1}.
\end{equation}
The values $U=1$ and $V=1$ do not occur: if one of them equals $1$, then \eqref{eq:fourier-relation} forces the other to equal $1$, and the middle coefficient $b$ vanishes.

Put
\[
 n_+:=(0,k_D)_4+(2,k_D)_4,\qquad n_-:=(1,k_D)_4+(3,k_D)_4,\qquad\varepsilon'_D:=(-1)^{r+s+1}.
\]
To see the indexing, put $x=U-1$.  Then $x+1=U\in C_{k_D}$, so $(j,k_D)_4$ counts the admissible $U$ with $x\in C_j$; the even values of $j$ give $\eta(U-1)=1$ and the odd values give $\eta(U-1)=-1$.  This is what $n_+$ and $n_-$ count.  Since $\eta(-\zeta)=(-1)^s$, equation \eqref{eq:D-ratio} turns the last condition into
\[
 \frac{\eta(V-1)}{\eta(U-1)}=\varepsilon'_D.
\]
Therefore
\begin{equation}\label{eq:Dcount}
 D_r=\frac{(n_++n_-)^2+\varepsilon'_D(n_+-n_-)^2}{2}.
\end{equation}
The pair $(U,V)$ determines $A_2/A_3$ by \eqref{eq:D-ratio}, and hence a unique projective coefficient class.  Conversely, each projective coefficient class of diagonal type gives one such pair, so no multiplicity is hidden in \eqref{eq:Dcount}.

\subsection{The edge type}
For edge type $+$ put
\[
 U:=A_0/A_1,\qquad V:=A_2/A_3,\qquad k_E\equiv r+2\pmod4.
\]
The conditions become
\[
 \chi(U)=\chi(V)=\zeta^{k_E},\qquad\eta(A_1/A_3)=-1,
\]
and \eqref{eq:fourier-relation} gives
\begin{equation}\label{eq:E-ratio}
 \frac{A_1}{A_3}=\frac{V+\zeta}{U+\zeta}.
\end{equation}
Suppose first that $U,V\neq-\zeta$, and define $i_E\equiv k_E-s\pmod4$ and
\begin{equation}\label{eq:mpm}
 m_+:=\sum_{\substack{j\in\mathbb Z/4\mathbb Z\\j\equiv s\ (2)}}(i_E,j)_4,\qquad
 m_-:=\sum_{\substack{j\in\mathbb Z/4\mathbb Z\\j\not\equiv s\ (2)}}(i_E,j)_4.
\end{equation}
For the indexing in \eqref{eq:mpm}, put $x=U/\zeta$.  Then $x\in C_{i_E}$ and $(U+\zeta)/\zeta=x+1$.  Since $\eta(\zeta)=(-1)^s$, one has $\eta(U+\zeta)=1$ exactly when the quartic class $j$ of $x+1$ satisfies $j\equiv s\pmod2$.  Thus $m_+$ and $m_-$ count the two signs of $\eta(U+\zeta)$.  The quotient in \eqref{eq:E-ratio} is a nonsquare precisely when these signs are opposite, so the classes with $U,V\neq-\zeta$ contribute $2m_+m_-$.

It remains to treat the case $U=-\zeta$.  If $U=-\zeta$, then \eqref{eq:fourier-relation} forces $V=-\zeta$.  Since $-\zeta=\gamma^{3s}$, its quartic class is $3s\bmod4$, so it lies in the required class $C_{k_E}$ exactly when
\begin{equation}\label{eq:exceptional-cong}
 r\equiv3s-2\pmod4.
\end{equation}
In that case put $\theta=A_1/A_3$.  The remaining condition is $\eta(\theta)=-1$, so $\theta$ is a free nonzero nonsquare.  The vector of values is proportional to $(-\zeta\theta,\theta,-\zeta,1)$, and its coefficient of $T$ is proportional to $1-\theta$, hence nonzero because a nonsquare is not equal to $1$.  Thus every such $\theta$ gives a unique projective coefficient class with $b\neq0$, and there are exactly $(q-1)/2$ of them.  With $\delta_r:=\mathbf 1_{r\equiv3s-2\ (4)}$ we obtain
\begin{equation}\label{eq:Ecount}
 E_r=2m_+m_-+\delta_r\frac{q-1}{2}.
\end{equation}
Combining \eqref{eq:D2E}, \eqref{eq:Dcount} and \eqref{eq:Ecount} gives
\begin{equation}\label{eq:core-count}
 N_{\mathrm{pre}}(r)=D_r+4m_+m_-+\delta_r(q-1).
\end{equation}

\subsection{Removing the binomials}\label{subsec:trinomial-count}
The count \eqref{eq:core-count} allows $a=0$ or $c=0$, and for these classes the exponent $r$ need not be the vanishing order.  The classes with $c=0$, $a\neq0$ contribute $B(r)$.  For the classes with $a=0$,
\[
 X^r(X^s+cX^{2s})=cX^{r+s}(X^s+c^{-1}),
\]
so for $a=0$, $c\neq0$ output scaling identifies the count with $B(r+s)$.  The two lines $a=0$ and $c=0$ meet at $a=c=0$, that is, at the monomial $X^{r+s}$.  The numbers $B(t)$ exclude a zero binomial coefficient, so this monomial belongs to neither of the two binomial counts and has to be added separately when it is a permutation monomial.  Put
\begin{equation}\label{eq:iota-r}
 \iota_r:=\mathbf 1_{\gcd(r+s,q-1)=1}.
\end{equation}
Then the union of the two lines $a=0$ and $c=0$ contains $B(r)+B(r+s)+\iota_r$ points counted by \eqref{eq:core-count}.  Since $\gcd(r,s)=1$ implies $\gcd(r+s,s)=1$, one has $\iota_r=\mathbf 1_{r+s\text{ odd}}$.

\begin{theorem}[Projective trinomial count]\label{thm:trinomial-count}
Let $q\equiv1\pmod4$, $s=(q-1)/4$, and $1\leqslant r\leqslant2s$ with $\gcd(r,s)=1$.  With $D_r$ as in \eqref{eq:Dcount}, $m_\pm$ as in \eqref{eq:mpm}, $\delta_r$ as in \eqref{eq:exceptional-cong} and $\iota_r$ as in \eqref{eq:iota-r},
\begin{equation}\label{eq:main-count}
 \mathfrak T(r)=D_r+4m_+m_-+\delta_r(q-1)-\mathcal B_{q,r}-\iota_r,
\end{equation}
where $\mathcal B_{q,r}$ is given by Theorem~\textup{\ref{thm:even-boundary}} if $q\equiv1\pmod8$ and by Theorem~\textup{\ref{thm:odd-boundary}} if $q\equiv5\pmod8$.
\end{theorem}

\begin{proof}
Equation \eqref{eq:core-count} gives the count of all classes with $b\neq0$.  Subtracting the points on the lines $a=0$ and $c=0$ gives \eqref{eq:main-count}, and the two evaluations of $\mathcal B_{q,r}$ were obtained in Sections~\ref{subsec:even-boundary} and~\ref{subsec:odd-boundary}.
\end{proof}

\begin{corollary}[Ordered coefficient triples]\label{cor:raw}
The number of ordered triples $(a,b,c)\in(\mathbb F_q^*)^3$ for which $X^r(a+bX^s+cX^{2s})$ permutes $\mathbb F_q$ is $(q-1)\,\mathfrak T(r)$.
\end{corollary}

\begin{proof}
Each projective class with $b\neq0$ has a unique representative with $b=1$, and multiplication by a nonzero scalar preserves the permutation property.
\end{proof}

\begin{proposition}[Index four]\label{prop:least-index}
Every trinomial counted by $\mathfrak T(r)$ has polynomial index exactly $4$.
\end{proposition}

\begin{proof}
Since $abc\neq0$, the exponent support is $\{r,r+s,r+2s\}$, so its difference gcd with $q-1=4s$ is $s$, and $\operatorname{ind}(F)=(q-1)/s=4$ \cite[\S15.1, pp.~319--320]{Wang2019}.  The bound $r\leqslant2s$ keeps the degree at most $q-1$.
\end{proof}

\begin{proof}[Proof of Theorem~\textup{\ref{thm:D}}]
Parts \ref{D:even} and \ref{D:odd} are Theorems~\ref{thm:even-boundary} and~\ref{thm:odd-boundary}, with $M_q$ and $\Lambda_q$ given by \eqref{eq:Mq} and \eqref{eq:Lambda-explicit}.  Part \ref{D:tri} is Theorem~\ref{thm:trinomial-count}, with $D_r$ and $m_\pm$ given by \eqref{eq:Dcount} and \eqref{eq:mpm}.  Part \ref{D:strata} was proved in Section~\ref{sec:index4}.
\end{proof}

\begin{remark}\label{rem:E128}
By Theorems~\ref{thm:even-boundary}, \ref{thm:odd-boundary} and~\ref{thm:trinomial-count}, the Frobenius trace $a_q(E_{128})$ enters the formula for $B(r)+B(r+s)$ with coefficient $-1/4$, and the formula for $\mathfrak T(r)$ with coefficient $1/4$.  In the sum of Theorem~\ref{thm:residue-collapse} these contributions cancel for each residue exponent separately: the Frobenius trace of the curve with $j$-invariant $128$, which has no complex multiplication over $\mathbb Q$, appears in the formulas for $B(r)+B(r+s)$ and $\mathfrak T(r)$ but not in their total.
\end{remark}


\begin{thebibliography}{99}

\bibitem{AGW2009}
A. Akbary, D. Ghioca and Q. Wang,
\emph{On permutation polynomials of prescribed shape},
Finite Fields Appl. 15 (2009), no. 2, 195--206.

\bibitem{AkbaryWang2007}
A. Akbary and Q. Wang,
\emph{On polynomials of the form $x^r f(x^{(q-1)/l})$},
Int. J. Math. Math. Sci. 2007 (2007), Art. ID 23408, 7 pp.

\bibitem{BiluLucaPizarro2016}
Y. Bilu, F. Luca and A. Pizarro-Madariaga,
\emph{Rational products of singular moduli},
J. Number Theory 158 (2016), 397--410.

\bibitem{BorsWang2022}
A. Bors and Q. Wang,
\emph{Generalized cyclotomic mappings: switching between polynomial, cyclotomic, and wreath product form},
Commun. Math. Res. 38 (2022), no. 2, 246--318.

\bibitem{DingZieve2023}
Z. Ding and M. E. Zieve,
\emph{Determination of a class of permutation quadrinomials},
Proc. Lond. Math. Soc. (3) 127 (2023), no. 2, 221--260.

\bibitem{DokchitserDokchitser2018}
T. Dokchitser and V. Dokchitser,
\emph{Quotients of hyperelliptic curves and \'etale cohomology},
Q. J. Math. 69 (2018), no. 2, 747--768.

\bibitem{FengWang2026}
X. Feng and Q. Wang,
\emph{Permutation polynomials of index $q+1$ over $\mathbb F_{q^2}$},
Finite Fields Appl. 112 (2026), Art. 102796.

\bibitem{Hou2025}
X.-d. Hou,
\emph{Number of equivalence classes of rational functions over finite fields},
Acta Arith. 218 (2025), no.~2, 97--136.

\bibitem{HouPallozziLavorante2023}
X.-d. Hou and V. Pallozzi Lavorante,
\emph{New results on permutation binomials of finite fields},
Finite Fields Appl. 88 (2023), Art.~102179.

\bibitem{Katz1988}
N. M. Katz, \emph{Gauss Sums, Kloosterman Sums, and Monodromy Groups},
Ann. of Math. Stud. 116, Princeton Univ. Press, 1988.

\bibitem{Katz2009Dwork}
N. M. Katz,
\emph{Another look at the Dwork family}, in Y. Tschinkel and Y. Zarhin (eds.),
\emph{Algebra, Arithmetic, and Geometry: In Honor of Yu. I. Manin. Vol. II},
Progr. Math. 270, Birkh\"auser Boston, 2009, pp. 89--126.

\bibitem{LeePark1997}
J. B. Lee and Y. H. Park,
\emph{Some permuting trinomials over finite fields},
Acta Math. Sci. 17 (1997), no. 3, 250--254.

\bibitem{LiGuanWu2025}
J. Li, H. Guan and D. Wu,
\emph{An enumeration problem for permutation polynomials},
Acta Math. Sinica (Chinese Ser.) 68 (2025), no. 5, 839--846.

\bibitem{Mullen1976}
G. L. Mullen,
\emph{Equivalence classes of polynomials over finite fields},
Acta Arith. 31 (1976), 113--123.

\bibitem{OliveiraBrochero2022}
J. A. Oliveira and F. E. Brochero Mart\'inez,
\emph{Permutation binomials over finite fields},
Discrete Math. 345 (2022), no.~3, Art.~112732.

\bibitem{ParkLee2001}
Y. H. Park and J. B. Lee,
\emph{Permutation polynomials and group permutation polynomials},
Bull. Aust. Math. Soc. 63 (2001), no. 1, 67--74.

\bibitem{Silverman2009}
J. H. Silverman,
\emph{The Arithmetic of Elliptic Curves}, 2nd ed., Grad. Texts in Math. 106, Springer, New York, 2009.

\bibitem{WanLidl1991}
D. Wan and R. Lidl,
\emph{Permutation polynomials of the form $x^r f(x^{(q-1)/d})$ and their group structure},
Monatsh. Math. 112 (1991), no. 2, 149--163.

\bibitem{Wang2019}
Q. Wang, \emph{Polynomials over finite fields: an index approach}, in K.-U. Schmidt and A. Winterhof (eds.),
\emph{Combinatorics and Finite Fields: Difference Sets, Polynomials, Pseudorandomness and Applications}, De Gruyter, Berlin/Boston, 2019, pp. 319--348.

\bibitem{WuYuanDingMa2017}
D. Wu, P. Yuan, C. Ding and Y. Ma,
\emph{Permutation trinomials over $\mathbb F_{2^m}$},
Finite Fields Appl. 46 (2017), 38--56.

\bibitem{YadavGuptaSinghYadav2026}
A. A. Yadav, I. Gupta, H. Singh and A. Yadav,
\emph{Several new classes of permutation polynomials over finite fields $\mathbb F_{5^{2m}}$},
Appl. Algebra Eng. Commun. Comput. 37 (2026), no. 4, 943--968.

\bibitem{YadavSinghGupta2025}
A. A. Yadav, H. Singh and I. Gupta,
\emph{A note on QM equivalence of known permutation polynomials},
Cryptogr. Commun. 17 (2025), no. 2, 525--540.

\end{thebibliography}
\end{document}